\documentclass[a4paper,11pt]{article}%
\usepackage{graphicx}
\usepackage{amsmath}
\usepackage{amsfonts}
\usepackage{amssymb}%
\newtheorem{theorem}{Theorem}

\newtheorem{corollary}[theorem]{Corollary}

\newtheorem{definition}{Definition}
\newtheorem{example}{Example}

\newtheorem{lemma}[theorem]{Lemma}
\newtheorem{notation}{Notation}

\newtheorem{proposition}[theorem]{Proposition}
\newtheorem{remark}{Remark}

\newenvironment{proof}[1][Proof]{\textbf{#1.} }{\ \rule{0.5em}{0.5em}}
\begin{document}

\title{Regular rapidly decreasing nonlinear generalized functions. Application to
microlocal regularity}
\author{Antoine Delcroix\\Equipe Analyse Alg\'{e}brique Non Lin\'{e}aire\\\textit{Laboratoire Analyse, Optimisation, Contr\^{o}le}\\Facult\'{e} des sciences - Universit\'{e} des Antilles et de la Guyane\\97159\ Pointe-\`{a}-Pitre Cedex Guadeloupe}
\maketitle

\begin{abstract}
We present new types of regularity for nonlinear generalized functions, based
on the notion of regular growth with respect to the regularizing parameter of
the Colombeau simplified model. This generalizes the notion of $\mathcal{G}%
^{\infty}$-regularity introduced by M.\ Oberguggenberger. \ A key point is
that these regularities can be characterized, for compactly supported
generalized functions, by a property of their Fourier transform. This opens
the door to microanalysis of singularities of generalized functions, with
respect to these regularities. We present a complete study of this topic,
including properties of the Fourier transform (exchange and regularity
theorems) and relationship with classical theory, via suitable results of embeddings.

\end{abstract}

\noindent\textbf{Mathematics Subject Classification (2000): 35A18, 35A27,
42B10, 46E10, 46F30}.\textbf{\smallskip}


\section{Introduction}

The various theories of nonlinear generalized functions are suitable
frameworks to set and solve differential or integral problems with irregular
operators or data. Even for linear problems, these theories are efficient to
overcome some limitations of the distributional framework. We follow in this
paper the theory introduced by J.-F.\ Colombeau \cite{Col1}, \cite{Col2},
\cite{GKOS}, \cite{NePiSc}. To be short, a special Colombeau type algebra is a
factor space $\mathcal{G}=\mathcal{X}/\mathcal{N}$ of moderate modulo
negligible nets.\ The moderateness (respectively the negligibility) of nets is
defined by their asymptotic behavior when a real parameter $\varepsilon$ tends
to $0$. \smallskip

A local and microlocal analysis of singularities of nonlinear generalized
functions has been developed during the last decade, based on the notion of
$\mathcal{G}^{\infty}$-regularity \cite{Ober1}. A generalized function is
$\mathcal{G}^{\infty}$-regular if it has uniform growth bounds, with respect
to the regularization parameter $\varepsilon$, for all derivatives. In fact,
this notion appears to be the exact generalization of the $\mathrm{C}^{\infty
}$-regularity for distributions, in the sense given by the result of
\cite{Ober1} which asserts that $\mathcal{G}^{\infty}\cap\mathcal{D}^{\prime}$
is equal to $\mathrm{C}^{\infty}$.\smallskip

In this paper we include $\mathcal{G}^{\infty}$ and $\mathcal{G}$ in a new
framework of spaces of $\mathcal{R}$-\textit{regular nonlinear generalized
functions}, in which the growth bounds are defined with the help of spaces
$\mathcal{R}$ of sequences satisfying natural conditions of stability. One
main property of those spaces is that the elements with compact support can be
characterized by a \textquotedblleft$\mathcal{R}$-property\textquotedblright%
\ of their Fourier transform.\ (Those Fourier Transforms belong to some
regular subspaces of spaces of \textit{rapidly decreasing generalized
functions} \cite{Garetto}, \cite{Radyno}.) Thus, the parallel is complete with
the $\mathrm{C}^{\infty}$-regularity of compactly supported distributions.
Moreover, from this characterization, we deduce that the microlocal behavior
of a generalized function with respect to a given $\mathcal{R}$-regularity is
completely similar to the one of a distribution with respect to the
$\mathrm{C}^{\infty}$-regularity.\ In particular, we can handle the
$\mathcal{R}$-wavefront of an element of $\mathcal{G}$ as the $\mathrm{C}%
^{\infty}$ one of a distribution. Finally, the $\mathcal{G}^{\infty}%
$-regularity for an element of $\mathcal{G}$ appears as a remarkable
particular case. \smallskip

With this new notion of $\mathcal{R}$-regularity, we enlarge the possibility
for the study of the propagation of singularities through differential and
pseudo differential operators, since we give a less restrictive framework than
the $\mathcal{G}^{\infty}$-regularity. This will be studied in forthcoming
papers.\smallskip

Let us also quote that spaces of $\mathcal{R}$-regular generalized functions
have been used in a problem of Schwartz kernel type theorem. More precisely,
we showed in \cite{DELGST} that some nets of linear maps (parametrized by
$\varepsilon\in\left(  0,1\right]  $), satisfying some growth conditions
similar to those introduced for $\mathcal{R}$-regular spaces, give rise to
linear maps between spaces of generalized functions. Moreover, those maps can
be represented by generalized integral kernel on some special $\mathcal{R}%
$-regular subspaces of $\mathcal{G}\left(  \Omega\right)  $ in which the
growth bounds are at most sublinear with respect to $l$. This reinforces the
interest of this framework.\smallskip

The paper is organized as follows. In section \ref{Glimsec2}, we introduce the
spaces of $\mathcal{R}$-regular generalized functions and we precise some
classical results about the embedding of $\mathcal{D}^{\prime}$ into these
spaces. Section \ref{Glimsec3} is devoted to the study of the space
$\mathcal{G}_{\mathcal{S}}$ of rapidly decreasing generalized functions.\ In
particular, we show that $\mathcal{O}_{C}^{^{\prime}}$, the space of rapidly
decreasing distributions, is embedded in $\mathcal{G}_{\mathcal{S}}$. Thus,
$\mathcal{G}_{\mathcal{S}}$ plays for $\mathcal{O}_{C}^{^{\prime}}$ the role
that $\mathcal{G}$ plays for $\mathcal{D}^{\prime}$. Section \ref{Glimsec4}
contains the material related to Fourier transform of elements of
$\mathcal{G}_{\mathcal{S}}$ and especially an exchange theorem which is, in
the context of $\mathcal{R}$-regularity, an analogon and a generalization of
the classical exchange theorem between $\mathcal{O}_{C}^{\prime}$ and
$\mathcal{O}_{M}$. Section \ref{Glimsec5} gives the above mentioned
characterization by Fourier transform of compactly supported $\mathcal{R}%
$-regular generalized functions whereas, in section \ref{Glimsec6}, we present
the $\mathcal{R}$-local and $\mathcal{R}$-microlocal analysis of generalized functions.

\section{The sheaf of Colombeau simplified algebras and related subsheaves
\label{Glimsec2}}

\subsection{Sheaves of regular generalized functions}

\begin{notation}
\label{GRegInfEg}For two sequences $\left(  N_{1},N_{2}\right)  \in\left(
\mathbb{R}_{+}^{\mathbb{N}}\right)  ^{2}$, we say that $N_{1}$ is smaller or
equal to $N_{2}$ and note $N_{1}\preccurlyeq N_{2}$ iff $\ \forall
n\in\mathbb{N}\ \ \ N_{1}\left(  n\right)  \leq N_{2}\left(  n\right)  $.
\end{notation}

\begin{definition}
\label{GlimDefReg}We say that a subspace $\mathcal{R}$ of $\mathbb{R}%
_{+}^{\mathbb{N}}$ is \emph{regular} if $\mathcal{R}$ is non empty
and\newline$\left(  i\right)  $~$\mathcal{R}$ is \textquotedblleft
overstable\textquotedblright\ by translation and by maximum$\vspace{-0.04in}$%
\begin{equation}
\forall N\in\mathcal{R},\ \ \forall\left(  k,k^{\prime}\right)  \in
\mathbb{N}^{2},\ \ \exists N^{\prime}\in\mathcal{R},\ \ \forall n\in
\mathbb{N}\ \ \ \ N\left(  n+k\right)  +k^{\prime}\leq N^{\prime}\left(
n\right)  ,\vspace{-0.04in} \label{GlimRegD1a}%
\end{equation}%
\begin{equation}
\forall N_{1}\in\mathcal{R},\ \ \forall N_{2}\in\mathcal{R},\ \ \exists
N\in\mathcal{R},\ \ \forall n\in\mathbb{N}\ \ \ \ \max\left(  N_{1}\left(
n\right)  ,N_{2}\left(  n\right)  \right)  \leq N\left(  n\right)  .
\label{GlimRegD1b}%
\end{equation}
$\left(  ii\right)  $~For all $N_{1}\ $and $N_{2}$ in $\mathcal{R}$, there
exists $N\in\mathcal{R}$ such that%
\begin{equation}
\forall\left(  l_{1},l_{2}\right)  \in\mathbb{N}^{2}\ \ \ \ \ N_{1}\left(
l_{1}\right)  +N_{2}\left(  l_{2}\right)  \leq N\left(  l_{1}+l_{2}\right)  .
\label{GlimRegD2}%
\end{equation}

\end{definition}

\begin{example}
\label{GlimExample0}~\newline$\left(  i\right)  $~The set $\mathcal{B}$ of
bounded sequences and the set $\mathcal{A}$ of affine sequences are regular
subsets of $\mathbb{R}_{+}^{\mathbb{N}}$, which is itself regular.\newline%
$\left(  ii\right)  $~The set $\mathcal{L}_{og}=\left\{  N\in\mathbb{R}%
_{+}^{\mathbb{N}}\,\left\vert \,\exists b\in\mathbb{R}_{+}\text{, }%
N:n\mapsto\ln n+b\right.  \right\}  $ is not regular ((\ref{GlimRegD2}) is not
satisfied), whereas $\mathcal{L}_{og}^{1}=\left\{  N\in\mathbb{R}%
_{+}^{\mathbb{N}}\,\left\vert \,\exists\left(  a,b\right)  \in\mathbb{R}%
_{+}^{2}\text{, }N:n\mapsto a\ln n+b\right.  \right\}  $ is regular.
((\ref{GlimRegD2}) comes, for example, from $\ln x+\ln y\leq2\ln(x+y)$, for
$x>0$ and $y>0$.)
\end{example}

Let $\Omega$ be an open subset of $\mathbb{R}^{d}$ ($d\in\mathbb{N}$) and
consider the algebra $\mathrm{C}^{\infty}\left(  \Omega\right)  $ of complex
valued smooth functions, endowed with its usual topology. This topology can be
described by the family of seminorms $\left(  p_{K,l}\right)  _{K\Subset
\Omega,l\in\mathbb{N}}$ defined by%
\[
p_{K,l}\left(  f\right)  =\sup_{x\in K,\left\vert \alpha\right\vert \leq
l}\left\vert \partial^{\alpha}f\left(  x\right)  \right\vert .
\]
For any regular subset $\mathcal{R}$ of $\mathbb{R}_{+}^{\mathbb{N}}$, we set
\begin{align*}
\mathcal{X}^{\mathcal{R}}\left(  \Omega\right)   &  =\left\{  \left(
f_{\varepsilon}\right)  _{\varepsilon}\in\mathrm{C}^{\infty}\left(
\Omega\right)  ^{\left(  0,1\right]  }\,\left\vert \,\forall K\Subset
\Omega,\;\exists N\in\mathcal{R},\;\forall l\in\mathbb{N}\;\;p_{K,l}\left(
f_{\varepsilon}\right)  =\mathrm{O}\left(  \varepsilon^{-N(l)}\right)
\;\mathrm{as}\;\varepsilon\rightarrow0\right.  \right\} \\
\mathcal{N}^{\mathcal{R}}\left(  \Omega\right)   &  =\left\{  \left(
f_{\varepsilon}\right)  _{\varepsilon}\in\mathrm{C}^{\infty}\left(
\Omega\right)  ^{\left(  0,1\right]  }\,\left\vert \,\forall K\Subset
\Omega,\;\forall m\in\mathcal{R},\ \forall l\in\mathbb{N}\;\;p_{K,l}\left(
f_{\varepsilon}\right)  =\mathrm{O}\left(  \varepsilon^{m(l)}\right)
\;\mathrm{as}\;\varepsilon\rightarrow0\right.  \right\}  .
\end{align*}

\begin{proposition}
\label{GlimPropFReg}~\newline$\left(  i\right)  $~For any regular subspace
$\mathcal{R}$ of $\mathbb{R}_{+}^{\mathbb{N}}$, the functor $\Omega
\rightarrow\mathcal{X}^{\mathcal{R}}\left(  \Omega\right)  $ defines a sheaf
of differential algebras over the ring
\[
\mathcal{X}\left(  \mathbb{C}\right)  =\left\{  \left(  r_{\varepsilon
}\right)  _{\varepsilon}\in\mathbb{C}^{\left(  0,1\right]  }\,\left\vert
\,\exists q\in\mathbb{N}\;\ \left\vert r_{\varepsilon}\right\vert
=\mathrm{O}\left(  \varepsilon^{-q}\right)  \;\mathrm{as}\;\varepsilon
\rightarrow0\right.  \right\}  .
\]
$\left(  ii\right)  $~The functor $\mathcal{N}^{\mathcal{R}}:\Omega
\rightarrow\mathcal{N}^{\mathcal{R}}\left(  \Omega\right)  $ defines a sheaf
of ideals of the sheaf $\mathcal{X}^{\mathcal{R}}\left(  \cdot\right)
$.\newline$\left(  iii\right)  $~For any regular subspaces $\mathcal{R}_{1}$
and $\mathcal{R}_{2}$ of $\mathbb{R}_{+}^{\mathbb{N}}$, with $\mathcal{R}%
_{1}\subset\mathcal{R}_{2}$, the sheaf $\mathcal{X}^{\mathcal{R}_{1}}\left(
\Omega\right)  $ is a subsheaf of the sheaf $\mathcal{X}^{\mathcal{R}_{2}%
}\left(  \Omega\right)  $.
\end{proposition}

\begin{proof}
We split the proof in two parts.\newline$\left(  a\right)  $%
~\textit{Algebraical properties.}-~Let us first show that for any open set
$\Omega\subset\mathbb{R}^{d}$, $\mathcal{X}^{\mathcal{R}}\left(
\Omega\right)  $ is a subalgebra of $\mathrm{C}^{\infty}\left(  \Omega\right)
^{\left(  0,1\right]  }$. Take $\left(  f_{\varepsilon}\right)  _{\varepsilon
}$ and $\left(  g_{\varepsilon}\right)  _{\varepsilon}$ in $\mathcal{X}%
^{\mathcal{R}}\left(  \Omega\right)  $ and $K\Subset\Omega$. There exist
$N_{f}\in\mathcal{R}$ and $N_{g}\in\mathcal{R}$ such that%
\[
\forall l\in\mathbb{N}\;\;\;p_{K,l}\left(  h_{\varepsilon}\right)
=\mathrm{O}\left(  \varepsilon^{-N_{h}\left(  l\right)  }\right)
\;\mathrm{as}\;\varepsilon\rightarrow0\text{, \ for }h_{\varepsilon
}=f_{\varepsilon},\,g_{\varepsilon}.
\]
We get immediately that $p_{K,l}\left(  f_{\varepsilon}+g_{\varepsilon
}\right)  =\mathrm{O}\left(  \varepsilon^{-\max\left(  N_{f}\left(  l\right)
,N_{g}\left(  l\right)  \right)  }\right)  $ as\ $\varepsilon\rightarrow0$,
with $\max\left(  N_{f},N_{g}\right)  \preccurlyeq N$ for some $N\in
\mathcal{R}$ according to (\ref{GlimRegD1b}). Then, $\left(  f_{\varepsilon
}+g_{\varepsilon}\right)  _{\varepsilon}$ belongs to $\mathcal{X}%
^{\mathcal{R}}\left(  \Omega\right)  .$

For $\left(  c_{\varepsilon}\right)  _{\varepsilon}\in\mathcal{X}\left(
\mathbb{C}\right)  $, there exists $q_{c}$ such that $\left\vert
c_{\varepsilon}\right\vert =\mathrm{O}\left(  \varepsilon^{-q_{c}}\right)  $
as\ $\varepsilon\rightarrow0$. Then $p_{K,l}\left(  c_{\varepsilon
}f_{\varepsilon}\right)  =\mathrm{O}\left(  \varepsilon^{-N_{f}\left(
l\right)  -q_{c}}\right)  $.\ From (\ref{GlimRegD1a}), there exists
$N\in\mathcal{R}$ such that $N_{f}+q_{c}\preccurlyeq N$.\ Thus, $\left(
c_{\varepsilon}f_{\varepsilon}\right)  _{\varepsilon}\in\mathcal{X}%
^{\mathcal{R}}\left(  \Omega\right)  $. It follows that $\mathcal{X}%
^{\mathcal{R}}\left(  \Omega\right)  $ is a submodule of $\mathrm{C}^{\infty
}\left(  \Omega\right)  ^{\left(  0,1\right]  }$ over $\mathcal{X}\left(
\mathbb{C}\right)  $.

Consider now $l\in\mathbb{N}$ and $\alpha\in\mathbb{N}^{d}$ with $\left\vert
\alpha\right\vert =l$. By Leibniz' formula, we have, for all $\varepsilon
\in\left(  0,1\right]  $ and $x\in K$,
\[
\left\vert \partial^{\alpha}\left(  f_{\varepsilon}g_{\varepsilon}\right)
\left(  x\right)  \right\vert =\sum_{\gamma\leq\alpha}C_{\alpha}^{\gamma
}\left\vert \partial^{\gamma}f_{\varepsilon}\,\left(  x\right)  \partial
^{\alpha-\gamma}g_{\varepsilon}\left(  x\right)  \right\vert \leq\sum
_{\gamma\leq\alpha}C_{\alpha}^{\gamma}p_{K,\left\vert \gamma\right\vert
}\left(  f_{\varepsilon}\right)  p_{K,\left\vert \alpha-\gamma\right\vert
}\left(  g_{\varepsilon}\right)  ,
\]
where $C_{\alpha}^{\gamma}$ is the generalized binomial coefficient. We have,
for all $\gamma\leq\alpha,$%
\[
p_{K,\left\vert \gamma\right\vert }\left(  f_{\varepsilon}\right)
p_{K,\left\vert \alpha-\gamma\right\vert }\left(  g_{\varepsilon}\right)
=\mathrm{O}\left(  \varepsilon^{-N_{f}\left(  \left\vert \gamma\right\vert
\right)  -Ng\left(  \left\vert \alpha-\gamma\right\vert \right)  }\right)
\ \mathrm{as}\;\varepsilon\rightarrow0.
\]
As $\gamma\leq\alpha$, we get $\left\vert \gamma\right\vert +\left\vert
\alpha-\gamma\right\vert =\left\vert \alpha\right\vert =l$. According to
(\ref{GlimRegD2}), there exists $N\in\mathcal{R}$ such that, for all $k$ and
$k^{\prime}\leq k$ in $\mathbb{N}$, $N_{f}\left(  k^{\prime}\right)
+N_{g}\left(  k-k^{\prime}\right)  \leq$ $N(k)$. Then
\[
\sup_{x\in K}\left\vert \partial^{\alpha}\left(  f_{\varepsilon}%
g_{\varepsilon}\right)  (x)\right\vert =\mathrm{O}\left(  \varepsilon
^{-N(l)}\right)  \ \mathrm{as}\;\varepsilon\rightarrow0.
\]
Thus $p_{K,l}\left(  f_{\varepsilon}\,g_{\varepsilon}\right)  =\mathrm{O}%
\left(  \varepsilon^{-N(l)}\right)  \ \mathrm{as}\;\varepsilon\rightarrow0$,
and $\left(  f_{\varepsilon}\,g_{\varepsilon}\right)  _{\varepsilon}%
\in\mathcal{X}^{\mathcal{R}}\left(  \Omega\right)  $.\smallskip

For the properties related to $\mathcal{N}^{\mathcal{R}}\left(  \Omega\right)
$, we have the following fundamental lemma:

\begin{lemma}
\label{GregLmnNRN}The set $\mathcal{N}^{\mathcal{R}}\left(  \Omega\right)  $
is equal to Colombeau's ideal%
\[
\mathcal{N}\left(  \Omega\right)  =\left\{  \left(  f_{\varepsilon}\right)
_{\varepsilon}\in\mathrm{C}^{\infty}\left(  \Omega\right)  ^{\left(
0,1\right]  }\,\left\vert \,\forall K\Subset\Omega,\;\forall l\in
\mathbb{N},\ \forall m\in\mathbb{N}\;\;p_{K,l}\left(  f_{\varepsilon}\right)
=\mathrm{O}\left(  \varepsilon^{m}\right)  \;\mathrm{as}\;\varepsilon
\rightarrow0\right.  \right\}  .
\]

\end{lemma}

Indeed, take $\left(  f_{\varepsilon}\right)  _{\varepsilon}\in\mathcal{N}%
^{\mathcal{R}}\left(  \Omega\right)  $. For any $K\Subset\Omega$%
,\ $l\in\mathbb{N}$ and $m\in\mathbb{N}$, choose $N\in\mathcal{R}$. According
to (\ref{GlimRegD1a}) there exists $N^{\prime}\in\mathcal{R}$ such that
$N+m\preccurlyeq N^{\prime}$. Thus, $p_{K,l}\left(  f_{\varepsilon}\right)
=\mathrm{O}\left(  \varepsilon^{N^{\prime}(l)}\right)  =\mathrm{O}\left(
\varepsilon^{m}\right)  $ as$\;\varepsilon\rightarrow0$ and $\left(
f_{\varepsilon}\right)  _{\varepsilon}\in\mathcal{N}\left(  \Omega\right)  $.
Conversely, given $\left(  f_{\varepsilon}\right)  _{\varepsilon}%
\in\mathcal{N}\left(  \Omega\right)  $ and $N\in\mathcal{R}$, we have
$p_{K,l}\left(  f_{\varepsilon}\right)  =\mathrm{O}\left(  \varepsilon
^{N(l)}\right)  $ as$\;\varepsilon\rightarrow0$, since this estimates holds
for all $m\in\mathbb{N}$.\smallskip

\noindent$\left(  b\right)  $~\textit{Sheaf properties}.-~The proof follows
the same lines as the one of Colombeau simplified algebras. (See for example
\cite{GKOS}, theorem 1.2.4.) First, the definition of restriction (by the mean
of the restriction of representatives) is straightforward as in Colombeau's
case. For the sheaf properties, we have to replace Colombeau's usual estimates
by $\mathcal{X}^{\mathcal{R}}$-estimates. But, at each place this happens, we
have only to consider a finite number of terms, by compactness properties.
Thus, the stability by maximum of $\mathcal{R}$ (property (\ref{GlimRegD1b}))
induces the result. Furthermore, lemma \ref{GregLmnNRN} shows that nothing
changes for estimates dealing with the ideal. Finally, point $\left(
iii\right)  $ of the proposition follows directly from the obvious inclusion
$\mathcal{X}^{\mathcal{R}_{1}}\left(  \Omega\right)  \subset\mathcal{X}%
^{\mathcal{R}_{2}}\left(  \Omega\right)  $.
\end{proof}

\begin{definition}
\label{GregularAlgebra}For any regular subset $\mathcal{R}$ of $\mathbb{R}%
_{+}^{\mathbb{N}}$, the sheaf of algebras
\[
\mathcal{G}^{\mathcal{R}}\left(  \,\cdot\,\right)  =\mathcal{X}^{\mathcal{R}%
}\left(  \,\cdot\,\right)  /\mathcal{N}^{\mathcal{R}}\left(  \,\cdot\,\right)
\]
is called the sheaf of $\mathcal{R}$-regular algebras of (nonlinear)
generalized functions.
\end{definition}

According to lemma \ref{GregLmnNRN}, we have $\mathcal{N}^{\mathcal{R}%
}=\mathcal{N}$. From now on, we shall write all algebras with Colombeau's
ideal $\mathcal{N}$. The sheaf $\mathcal{G}^{\mathcal{R}}\left(  \cdot\right)
$ turns out to be a sheaf of differential algebras and a sheaf of modules over
the factor ring $\overline{\mathbb{C}}=\mathcal{X}\left(  \mathbb{C}\right)
/\mathcal{N}\left(  \mathbb{C}\right)  $ with
\[
\mathcal{N}\left(  \mathbb{C}\right)  =\left\{  \left(  r_{\varepsilon
}\right)  _{\varepsilon}\in\mathbb{C}^{\left(  0,1\right]  }\,\left\vert
\,\forall p\in\mathbb{N}\;\ \left\vert r_{\varepsilon}\right\vert
=\mathrm{O}\left(  \varepsilon^{p}\right)  \;\mathrm{as}\;\varepsilon
\rightarrow0\right.  \right\}  .
\]

\begin{example}
\label{GlimExample1}Taking $\mathcal{R}=\mathbb{R}_{+}^{\mathbb{N}}$, we
recover the sheaf of \emph{Colombeau simplified }or\emph{ special algebras}.
\end{example}

\begin{notation}
\label{NotaGclass}In the sequel, we shall write $\mathcal{G}\left(
\Omega\right)  $ (resp. $\mathcal{X}_{M}\left(  \Omega\right)  $) instead of
$\mathcal{G}^{\mathbb{R}_{+}^{\mathbb{N}}}\left(  \Omega\right)  $ (resp.
$\mathcal{X}^{\mathbb{R}_{+}^{\mathbb{N}}}\left(  \Omega\right)  $). For
$\left(  f_{\varepsilon}\right)  _{\varepsilon}$ in $\mathcal{X}_{M}\left(
\Omega\right)  $ or $\mathcal{X}^{\mathcal{R}}\left(  \Omega\right)  $,
$\left[  \left(  f_{\varepsilon}\right)  _{\varepsilon}\right]  $ will be its
class in $\mathcal{G}\left(  \Omega\right)  $ or in $\mathcal{G}^{\mathcal{R}%
}\left(  \Omega\right)  $, since these classes are obtained modulo the same
ideal. (We consider $\mathcal{G}^{\mathcal{R}}\left(  \Omega\right)  $ as a
subspace of $\mathcal{G}\left(  \Omega\right)  $.)
\end{notation}

\begin{example}
\label{GlimExample2}Taking $\mathcal{R}=\mathcal{B}$, introduced in example
\ref{GlimExample0}, we obtain the sheaf of $\mathcal{G}^{\infty}$-generalized
functions \cite{Ober1}.
\end{example}

\begin{example}
\label{Glimexample3}Take $a$ in $\left[  0,+\infty\right]  $ and set
\[
\mathcal{R}_{0}=\left\{  N\in\mathbb{R}_{+}^{\mathbb{N}}\,\left\vert
\,\lim_{l\rightarrow+\infty}\left(  N(l)/l\right)  =0\right.  \right\}
\ ;\text{\ For }a>0:\mathcal{R}_{a}=\left\{  N\in\mathbb{R}_{+}^{\mathbb{N}%
}\,\left\vert \,\underset{l\rightarrow+\infty}{\lim\!\sup}\left(
N(l)/l\right)  <a\right.  \right\}  \text{.}%
\]
For any $a$ in $\left[  0,+\infty\right]  $, $\mathcal{R}_{a}$ is a regular
subset of $\mathbb{R}_{+}^{\mathbb{N}}$. The corresponding sheaves
$\mathcal{G}^{\mathcal{R}_{a}}\left(  \cdot\right)  $ are the sheaves of
algebras of generalized functions with slow growth introduced in
\cite{DELGST}. \ Note that, for $a$ in $\left(  0,+\infty\right]  $, a
sequence $N$ is in $\mathcal{R}_{a}$ iff there exists $\left(  a^{\prime
},b\right)  \in\left(  \mathbb{R}^{+}\right)  ^{2}$ with $a^{\prime}<a$ such
that $N(l)\leq a^{\prime}l+b$. The growth of the sequence $N$ is at most linear.
\end{example}

\begin{remark}
\label{GRegInclusion}We have the obvious sheaves inclusions$\ \mathcal{X}%
^{\infty}\left(  \,\cdot\,\right)  \subset\mathcal{X}^{\mathcal{R}_{a}}\left(
\,\cdot\,\right)  \subset\mathcal{X}\left(  \,\cdot\,\right)  .$
\end{remark}

For a given regular subspace $\mathcal{R}$ of $\mathbb{R}_{+}^{\mathbb{N}}$,
the notion of support of a section $f\in\mathcal{G}^{\mathcal{R}}\left(
\Omega\right)  $ ($\Omega$ open subset of $\mathbb{R}^{d}$) makes sense since
$\mathcal{G}^{\mathcal{R}}\left(  \cdot\right)  $ is a sheaf. The following
definition will be sufficient for this paper.

\begin{definition}
\label{GRegSuppKomp}The support of a generalized function $f\in\mathcal{G}%
^{\mathcal{R}}\left(  \Omega\right)  $ is the complement in $\Omega$ of the
largest open subset of $\Omega$ where $f$ is null.
\end{definition}

\begin{notation}
\label{NotGRegSuppKomp}We denote by $\mathcal{G}_{C}^{\mathcal{R}}\left(
\Omega\right)  $ the subset of $\mathcal{G}^{\mathcal{R}}\left(
\Omega\right)  $ of elements with compact support.
\end{notation}

\begin{lemma}
\label{LmnSuppK}Every $f\in\mathcal{G}_{C}^{\mathcal{R}}$ has a representative
$\left(  f_{\varepsilon}\right)  _{\varepsilon}$, such that each
$f_{\varepsilon}$ has the same compact support.
\end{lemma}

We shall not prove this lemma here, since lemma \ref{GregEmbCRSRu} below gives
the main ideas of the proof.

\subsection{Some embeddings}

For any regular subspace $\mathcal{R}$ of $\mathbb{R}_{+}^{\mathbb{N}}$ and
any $\Omega$ open subset of $\mathbb{R}^{d}$, $\mathrm{C}^{\infty}\left(
\Omega\right)  $ is embedded into $\mathcal{G}^{\mathcal{R}}\left(
\Omega\right)  $ by the canonical embedding
\[
\sigma:\mathrm{C}^{\infty}\left(  \Omega\right)  \rightarrow\mathcal{G}%
^{\mathcal{R}}\left(  \Omega\right)  \;\;\ f\rightarrow\left[  \left(
f_{\varepsilon}\right)  _{\varepsilon}\right]  \;\mathrm{with\;}%
f_{\varepsilon}=f\text{ for all }\varepsilon\in\left(  0,1\right]  .
\]
We refine here the well known result about the embedding of $\mathcal{D}%
^{\prime}\left(  \Omega\right)  $ into $\mathcal{G}\left(  \Omega\right)
$.\ Following the ideas of \cite{NePiSc}, we consider $\rho\in\mathcal{S}%
\left(  \mathbb{R}^{d}\right)  $ such that
\[%
{\textstyle\int}
\rho\left(  x\right)  \,\mathrm{d}x=1,\ \ \ \ \ \ \ \ \ \ \
{\textstyle\int}
x^{m}\rho\left(  x\right)  \,\mathrm{d}x=0\;\mathrm{for\;all\;}m\in
\mathbb{N}^{d}\,\backslash\left\{  0\right\}  .
\]
(Such a map can be chosen as the Fourier transform of a function of
$\mathcal{D}\left(  \mathbb{R}^{d}\right)  $ equal to $1$ on a neighborhood of
$0$.)\ We now choose $\chi\in\mathcal{D}\left(  \mathbb{R}^{d}\right)  $ such
that $0\leq\chi\leq1$,$\;\chi\equiv1$ on $\overline{B(0,1)}$ and $\chi\equiv0$
on $\mathbb{R}^{d}\backslash B(0,2)$.\ We define
\[
\forall\varepsilon\in\left(  0,1\right]  ,\;\;\forall x\in\mathbb{R}%
^{d}\;\;\;\theta_{\varepsilon}\left(  x\right)  =\frac{1}{\varepsilon^{d}%
}\,\rho\left(  \frac{x}{\varepsilon}\right)  \chi\left(  \left\vert
\ln\varepsilon\right\vert x\right)  \,\text{.}%
\]
Finally, consider $\left(  \kappa_{\varepsilon}\right)  _{\varepsilon}%
\in\left(  \mathcal{D}\left(  \mathbb{R}^{d}\right)  \right)  ^{\left(
0,1\right]  }$ such that
\[
\forall\varepsilon\in\left(  0,1\right)  ,\;\;0\leq\kappa_{\varepsilon}%
\leq1\;\;\;\kappa_{\varepsilon}\equiv1\text{ on }\left\{  x\in\Omega
\,\left\vert \,d(x,\mathbb{R}^{d}\,\text{%
$\backslash$%
}\Omega)\geq\varepsilon\text{ and }d(x,0)\leq1/\varepsilon\right.  \right\}
.
\]
With these ingredients, the map
\begin{equation}
\iota:\mathcal{D}^{\prime}\left(  \Omega\right)  \rightarrow\mathcal{G}\left(
\Omega\right)  \;\;\;T\mapsto\left(  \kappa_{\varepsilon}T\ast\theta
_{\varepsilon}\right)  _{\varepsilon}+\mathcal{N}\left(  \Omega\right)
\label{GrDefMollyDP}%
\end{equation}
is an embedding of $\mathcal{D}^{\prime}\left(  \Omega\right)  $ into
$\mathcal{G}\left(  \Omega\right)  $ such that $\iota_{\left\vert
\mathrm{C}^{\infty}\left(  \Omega\right)  \right.  }=\sigma$.

The proof is mainly based on the following property of $\left(  \theta
_{\varepsilon}\right)  _{\varepsilon}$:
\begin{equation}%
{\textstyle\int}
\theta_{\varepsilon}\left(  x\right)  \,\mathrm{d}x=1+\mathrm{O}\left(
\varepsilon^{k}\right)  \;\text{as}\;\varepsilon\rightarrow0,\ \ \ \forall
m\in\mathbb{N}^{d}\,\backslash\left\{  0\right\}  \;\
{\textstyle\int}
x^{m}\theta_{\varepsilon}\left(  x\right)  \,\mathrm{d}x=\mathrm{O}\left(
\varepsilon^{k}\right)  \;\text{as}\;\varepsilon\rightarrow0.
\label{Gregmollif}%
\end{equation}

Set
\begin{equation}
\mathcal{R}_{1}=\left\{  N\in\mathbb{R}_{+}^{\mathbb{N}}\,\left\vert \,\exists
b\in\mathbb{R}_{+},\forall l\in\mathbb{R}\ \ N(l)\leq l+b\right.  \right\}  .
\label{GregDefR1}%
\end{equation}
One can verify that the set $\mathcal{R}$ is regular and we set $\mathcal{G}%
^{\left(  1\right)  }\left(  \cdot\right)  =$ $\mathcal{G}^{\mathcal{R}_{1}%
}\left(  \cdot\right)  $.

\begin{proposition}
\label{LmnEmbedG1}The image of $\mathcal{D}^{\prime}\left(  \Omega\right)  $
by the embedding, defined by (\ref{GrDefMollyDP}), is included in
$\mathcal{G}^{\left(  1\right)  }\left(  \Omega\right)  $.
\end{proposition}

\begin{proof}
The proof is a refinement of the classical proof (see \cite{DelEmb}), and we
shall focus on the estimation of the growth of $\left(  \kappa_{\varepsilon
}T\ast\theta_{\varepsilon}\right)  _{\varepsilon}$ which is the main novelty.
We shall do it for the case $\Omega=\mathbb{R}^{d}$, for which the additional
cutoff by $\left(  \kappa_{\varepsilon}\right)  _{\varepsilon}$ is not needed.
Then, for a given $T\in\mathcal{D}^{\prime}\left(  \mathbb{R}^{d}\right)  $,
we have $\iota\left(  T\right)  =\left(  T\ast\theta_{\varepsilon}\right)
_{\varepsilon}+\mathcal{N}\left(  \mathbb{R}^{d}\right)  $, with $\left(
T\ast\theta_{\varepsilon}\right)  _{\varepsilon}=\left\langle T,\theta
_{\varepsilon}\left(  y-\cdot\right)  \right\rangle $.

Let us fix a compact set $K$ and consider $W$ an open subset of $\mathbb{R}%
^{d}$ such that $K\subset W\subset\overline{W}\Subset\mathbb{R}^{d}$. With the
above definitions, the function $x\mapsto\theta_{\varepsilon}\left(
y-x\right)  $ belongs to $\mathcal{D}\left(  W\right)  $ for all $y\in K$ and
$\varepsilon$ small enough, since the support of $\theta_{\varepsilon}$
shrinks to $\left\{  0\right\}  $ when $\varepsilon$ tends to $0$. Therefore,
for $\beta\in\mathbb{N}^{d}$ and $\varepsilon$ small enough, we have
\begin{align*}
\left.  \forall y\in K\ \ \ \ \partial^{\beta}\left(  T\ast\theta
_{\varepsilon}\right)  (y)\right.   &  =\left\langle T,\partial^{\beta
}\left\{  x\mapsto\theta_{\varepsilon}\left(  y-x\right)  \right\}
\right\rangle \\
&  =\left\langle T_{\left\vert W\right.  },\partial^{\beta}\left\{
x\mapsto\theta_{\varepsilon}\left(  y-x\right)  \right\}  \right\rangle \\
&  =\left(  -1\right)  ^{\left\vert \beta\right\vert }\left\langle
T_{\left\vert W\right.  },\left\{  x\mapsto\left(  \partial^{\beta}%
\theta_{\varepsilon}\right)  \left(  y-x\right)  \right\}  \right\rangle .
\end{align*}
By using the local structure of distributions \cite{Schwartz1}, we can write
$T_{\left\vert W\right.  }=\partial_{x}^{\alpha}f$ where $f$ is a compactly
supported continuous function having its support included in $W$. It follows
\begin{align*}
\left.  \forall y\in K\ \ \ \ \partial^{\beta}\left(  T\ast\theta
_{\varepsilon}\right)  (y)\right.   &  =\left(  -1\right)  ^{\left\vert
\beta\right\vert }\left\langle \partial_{x}^{\alpha}f,\left(  \partial^{\beta
}\theta_{\varepsilon}\right)  \left(  y-\cdot\right)  \right\rangle \\
&  =\left(  -1\right)  ^{\left\vert \alpha\right\vert +\left\vert
\beta\right\vert }\left\langle f,\left(  \partial^{\alpha+\beta}%
\theta_{\varepsilon}\right)  \left(  y-\cdot\right)  \right\rangle \\
&  =\left(  -1\right)  ^{\left\vert \alpha\right\vert +\left\vert
\beta\right\vert }\int_{W}f\left(  x\right)  \partial^{\alpha+\beta}%
\theta_{\varepsilon}\left(  y-x\right)  \,\mathrm{d}x.
\end{align*}
Using the definition of $\left(  \theta_{\varepsilon}\right)  _{\varepsilon}$,
we get
\begin{align*}
\left.  \forall\xi\in\mathbb{R}\;\;\;\partial^{\alpha+\beta}\theta
_{\varepsilon}\left(  \xi\right)  \right.   &  =\sum_{\gamma\leq\alpha+\beta
}\mathrm{C}_{\alpha+\beta}^{\gamma}\,\partial^{\gamma}\rho_{\varepsilon
}\left(  \xi\right)  \,\partial^{\alpha+\beta-\gamma}\left(  \chi\left(
\xi\left\vert \ln\varepsilon\right\vert \right)  \right)  \text{ (with }%
\rho_{\varepsilon}\left(  \cdot\right)  =\frac{1}{\varepsilon^{d}}\rho\left(
\frac{\cdot}{\varepsilon}\right)  \text{)}\\
&  =\sum_{\gamma\leq\alpha+\beta}\mathrm{C}_{\alpha+\beta}^{\gamma
}\,\varepsilon^{-d-\left\vert \gamma\right\vert }\left\vert \ln\varepsilon
\right\vert ^{\left\vert \alpha+\beta-\gamma\right\vert }\,\left(
\partial^{\gamma}\rho\right)  \left(  \frac{\xi}{\varepsilon}\right)  \left(
\partial^{\alpha+\beta-\gamma}\chi\right)  \left(  \xi\left\vert
\ln\varepsilon\right\vert \right)  ,
\end{align*}
with $\left\vert \alpha+\beta-\gamma\right\vert =\left\vert \alpha\right\vert
+\left\vert \beta\right\vert -\left\vert \gamma\right\vert $ since $\gamma
\leq\alpha+\beta$. For all $\gamma\leq\alpha+\beta$, we have
\[
\varepsilon^{-d-\left\vert \gamma\right\vert }\left\vert \ln\varepsilon
\right\vert ^{\left\vert \alpha+\beta-\gamma\right\vert }=\mathrm{O}\left(
\varepsilon^{-d-1-\left\vert \alpha\right\vert -\left\vert \beta\right\vert
}\right)  \text{ as }\varepsilon\rightarrow0.
\]
As $\rho$ and $\chi$ are bounded, as well as their derivatives, there exists
$C_{1}>0$, depending on $\left(  \alpha,\beta\right)  $, such that
\[
\forall\xi\in\mathbb{R}\;\;\;\left\vert \partial^{\alpha+\beta}\left(
\theta_{\varepsilon}\left(  \xi\right)  \right)  \right\vert \leq
C_{1}\,\varepsilon^{-d-1-\left\vert \alpha\right\vert -\left\vert
\beta\right\vert }.
\]
We get
\[
\forall y\in K\;\;\;\left\vert \partial^{\beta}\left(  T\ast\theta
_{\varepsilon}\right)  (y)\right\vert \leq C_{1}\,\sup_{\xi\in\overline{W}%
}\left\vert f\left(  \xi\right)  \right\vert \,\mathrm{vol}\left(
\overline{W}\right)  \,\varepsilon^{-d-1-\left\vert \alpha\right\vert
-\left\vert \beta\right\vert }.
\]
From this last inequality, we deduce that for any $l\in\mathbb{N}$, there
exists some constant $C_{2}>0$, depending on $\left(  l,K\right)  $, such
that
\[
p_{K,l}\left(  \left(  T\ast\theta_{\varepsilon}\right)  \right)  \leq
C_{2}\varepsilon^{-d-1-\left\vert \alpha\right\vert -l}=C_{2}\,\varepsilon
^{-N(l)}\text{ with }N(l)=l+d+1+\left\vert \alpha\right\vert ,
\]
which shows that $\left(  T\ast\theta_{\varepsilon}\right)  _{\varepsilon}%
\in\mathcal{X}^{\mathcal{R}_{1}}\left(  \mathbb{R}^{d}\right)  $%
.\medskip\smallskip
\end{proof}

We can summarize these results in the following commutative diagram:%

\begin{equation}%
\begin{array}
[c]{ccccc}%
\mathrm{C}^{\infty}\left(  \Omega\right)  & \longrightarrow & \mathcal{D}%
^{\prime}\left(  \Omega\right)  &  & \\
\downarrow\sigma &  & \downarrow\iota &  & \\
\mathcal{G}^{\infty}\left(  \Omega\right)  & \longrightarrow & \mathcal{G}%
^{\left(  1\right)  }\left(  \Omega\right)  & \longrightarrow & \mathcal{G}%
\left(  \Omega\right)  \,.
\end{array}
\label{ScomDiagDP}%
\end{equation}

\section{Rapidly decreasing generalized functions \label{Glimsec3}}

\subsection{Definition and first properties}

Spaces of rapidly decreasing generalized functions have been introduced in the
literature (\cite{Garetto}, \cite{Radyno}, \cite{Scarpa1}), notably in view of
the definition of the Fourier transform in convenient spaces of nonlinear
generalized functions. We give here a more complete description of this type
of space in the framework of $\mathcal{R}$-regular spaces.

\begin{definition}
\label{GlimDefReg2D}We say that a subspace $\mathcal{R}^{\prime}$ of the space
$\mathbb{R}_{+}^{\mathbb{N}^{2}}$of maps from $\mathbb{N}^{2}$ to
$\mathbb{R}_{+}$ is \emph{regular} if\newline$\left(  i\right)  $%
~$\mathcal{R}^{\prime}$ is \textquotedblleft overstable\textquotedblright\ by
translation and by maximum%
\begin{equation}
\forall N\in\mathcal{R}^{\prime},\ \forall\left(  k,k^{\prime},k^{\prime
\prime}\right)  \in\mathbb{N}^{3},\ \exists N^{\prime}\in\mathcal{R}^{\prime
},\ \forall\left(  q,l\right)  \in\mathbb{N}^{2}\ \ \ N\left(  q+k,l+k^{\prime
}\right)  +k^{\prime\prime}\leq N^{\prime}\left(  q,l\right)  ,\vspace
{-0.04in} \label{GlimRegD2D1a}%
\end{equation}%
\begin{equation}
\forall N_{1}\in\mathcal{R}^{\prime},\ \forall N_{2}\in\mathcal{R}^{\prime
},\ \exists N\in\mathcal{R}^{\prime},\ \forall\left(  q,l\right)
\in\mathbb{N}^{2}\ \ \ \ \max\left(  N_{1}(q,l),N_{2}(q,l)\right)  \leq
N(q,l). \label{GlimRegD2D1b}%
\end{equation}
$\left(  ii\right)  $~For any $N_{1}\ $and $N_{2}$ in $\mathcal{R}^{\prime}$,
there exists $N\in\mathcal{R}^{\prime}$ such that%
\begin{equation}
\forall\left(  q_{1},q_{2},l_{1},l_{2}\right)  \in\mathbb{N}^{4}%
\ \ \ \ \ N_{1}\left(  q_{1},l_{1}\right)  +N_{2}\left(  q_{2},l_{2}\right)
\leq N\left(  q_{1}+q_{2},l_{1}+l_{2}\right)  . \label{GlimRegD2D2}%
\end{equation}

\end{definition}

\begin{example}
\label{GlimExample0D2}~\newline$\left(  i\right)  $~The set $\mathcal{B}%
^{\prime}$ of bounded maps from $\mathbb{N}^{2}$ to $\mathbb{R}_{+}$ is a
regular subset of $\mathbb{R}_{+}^{\mathbb{N}^{2}}$.\newline$\left(
ii\right)  $~The set $\mathbb{R}_{+}^{\mathbb{N}^{2}}$ of all maps from
$\mathbb{N}^{2}$ to $\mathbb{R}_{+}$ is a regular set.
\end{example}

We consider $\Omega$ an open subset of $\mathbb{R}^{d}$ and the space
$\mathcal{S}\left(  \Omega\right)  $ of rapidly decreasing functions defined
on $\Omega$, endowed with the family of seminorms $\mathcal{Q}\left(
\Omega\right)  =\left(  \mu_{q,l}\right)  _{\left(  q,l\right)  \in
\mathbb{N}^{2}}$ defined by
\[
\mu_{q,l}\left(  f\right)  =\sup_{x\in\Omega,\left\vert \alpha\right\vert \leq
l}\left(  1+\left\vert x\right\vert \right)  ^{q}\left\vert \partial^{\alpha
}f\left(  x\right)  \right\vert .
\]
Let $\mathcal{R}^{\prime}$ be a regular subset of $\mathbb{R}_{+}%
^{\mathbb{N}^{2}}$ and set
\begin{align*}
\mathcal{X}_{\mathcal{S}}^{\mathcal{R}^{\prime}}\left(  \Omega\right)   &
=\left\{  \left(  f_{\varepsilon}\right)  _{\varepsilon}\in\mathcal{S}\left(
\Omega\right)  ^{\left(  0,1\right]  }\,\left\vert \,\exists N\in
\mathcal{R}^{\prime},\;\forall\left(  q,l\right)  \in\mathbb{N}^{2}%
\;\;\mu_{q,l}\left(  f_{\varepsilon}\right)  =\mathrm{O}\left(  \varepsilon
^{-N\left(  q,l\right)  }\right)  \;\mathrm{as}\;\varepsilon\rightarrow
0\right.  \right\}  ,\\
\mathcal{N}_{\mathcal{S}}^{\mathcal{R}^{\prime}}\left(  \Omega\right)   &
=\left\{  \left(  f_{\varepsilon}\right)  _{\varepsilon}\in\mathcal{S}\left(
\Omega\right)  ^{\left(  0,1\right]  }\,\left\vert \,\forall m\in
\mathcal{R}^{\prime},\ \forall\left(  q,l\right)  \in\mathbb{N}^{2}%
\;\;\mu_{q,l}\left(  f_{\varepsilon}\right)  =\mathrm{O}\left(  \varepsilon
^{m(q,l)}\right)  \;\mathrm{as}\;\varepsilon\rightarrow0\right.  \right\}  .
\end{align*}
As for lemma \ref{GregLmnNRN}, we have $\mathcal{N}_{\mathcal{S}}%
^{\mathcal{R}^{\prime}}\left(  \Omega\right)  =\mathcal{N}_{\mathcal{S}%
}\left(  \Omega\right)  $, with%
\[
\mathcal{N}_{\mathcal{S}}\left(  \Omega\right)  =\left\{  \left(
f_{\varepsilon}\right)  _{\varepsilon}\in\mathcal{S}\left(  \Omega\right)
^{\left(  0,1\right]  }\,\left\vert \,\forall\left(  q,l\right)  \in
\mathbb{N}^{2},\;\forall m\in\mathbb{N}\ \ \;\mu_{q,l}\left(  f_{\varepsilon
}\right)  =\mathrm{O}\left(  \varepsilon^{m}\right)  \;\mathrm{as}%
\;\varepsilon\rightarrow0\right.  \right\}  .
\]

\begin{proposition}
\label{GlimPropFRegS}~\newline$\left(  i\right)  $~For any regular subspace
$\mathcal{R}^{\prime}$ of $\mathbb{R}_{+}^{\mathbb{N}^{2}}$, the functor
$\Omega\rightarrow\mathcal{X}_{\mathcal{S}}^{\mathcal{R}^{\prime}}\left(
\Omega\right)  $ defines a presheaf (it allows restrictions) of differential
algebras over the ring $\mathcal{X}\left(  \mathbb{C}\right)  $.\newline%
$\left(  ii\right)  $~The functor $\mathcal{N}_{\mathcal{S}}:\Omega
\rightarrow\mathcal{N}_{\mathcal{S}}\left(  \Omega\right)  $ defines a
presheaf of ideals of the presheaf $\mathcal{X}_{\mathcal{S}}^{\mathcal{R}%
^{\prime}}\left(  \cdot\right)  $.\newline$\left(  iii\right)  $~For any
regular subspaces $\mathcal{R}_{1}^{\prime}$ and $\mathcal{R}_{2}^{\prime}$ of
$\mathbb{R}_{+}^{\mathbb{N}^{2}}$, with $\mathcal{R}_{1}^{\prime}%
\subset\mathcal{R}_{2}^{\prime}$, the presheaf $\mathcal{X}_{\mathcal{S}%
}^{\mathcal{R}_{1}^{\prime}}\left(  \Omega\right)  $ is a subpresheaf of the
presheaf $\mathcal{X}_{\mathcal{S}}^{\mathcal{R}_{2}^{\prime}}\left(
\Omega\right)  $.
\end{proposition}

\begin{proof}
$~$\newline$\left(  a\right)  $~\textit{Algebraical properties}.-~Let us first
prove that $\mathcal{X}_{\mathcal{S}}^{\mathcal{R}^{\prime}}\left(
\Omega\right)  $ is a subalgebra of $\mathcal{S}\left(  \Omega\right)
^{\left(  0,1\right]  }$. The proof that $\mathcal{X}_{\mathcal{S}%
}^{\mathcal{R}^{\prime}}\left(  \Omega\right)  $ is a sublinear space of
$\mathrm{C}^{\infty}\left(  \Omega\right)  ^{\left(  0,1\right]  }$ goes along
the same lines as in proposition \ref{GlimPropFReg}. For the product, take
$\left(  f_{\varepsilon}\right)  _{\varepsilon}$ and $\left(  g_{\varepsilon
}\right)  _{\varepsilon}$ in $\mathcal{X}_{\mathcal{S}}^{\mathcal{R}^{\prime}%
}\left(  \Omega\right)  $. According to the definitions, there exist $N_{f}%
\in\mathcal{R}^{\prime}$ and $N_{g}\in\mathcal{R}^{\prime}$ such that%
\[
\forall\left(  q,l\right)  \in\mathbb{N}^{2}\;\;\;\mu_{q,l}\left(
h_{\varepsilon}\right)  =\mathrm{O}\left(  \varepsilon^{-N_{h}\left(
q,l\right)  }\right)  \;\mathrm{as}\;\varepsilon\rightarrow0\text{, for
}h_{\varepsilon}=f_{\varepsilon},\,g_{\varepsilon}.
\]
Consider $\left(  q,l\right)  \in\mathbb{N}^{2}$ and $\alpha\in\mathbb{N}^{d}$
with $\left\vert \alpha\right\vert =l$. By Leibniz' formula, we have
\[
\forall\varepsilon\in\left(  0,1\right]  \text{\ \ \ \ \ }\partial^{\alpha
}\left(  f_{\varepsilon}g_{\varepsilon}\right)  =\sum_{\gamma\leq\alpha
}C_{\alpha}^{\gamma}\,\partial^{\gamma}f_{\varepsilon}\,\partial
^{\alpha-\gamma}g_{\varepsilon}.
\]
Thus%
\[
\sup_{x\in\Omega}\left(  1+\left\vert x\right\vert \right)  ^{q}\left\vert
\partial^{\alpha}\left(  f_{\varepsilon}g_{\varepsilon}\right)  (x)\right\vert
\leq\sum_{\gamma\leq\alpha}C_{\alpha}^{\gamma}\,\mu_{q,\left\vert
\gamma\right\vert }\left(  f_{\varepsilon}\right)  \,\mu_{0,\left\vert
\alpha-\gamma\right\vert }\left(  f_{\varepsilon}\right)  ,
\]
with, for all $\gamma\leq\alpha$, $\mu_{q,\left\vert \gamma\right\vert
}\left(  f_{\varepsilon}\right)  \mu_{0,\left\vert \alpha-\gamma\right\vert
}\left(  f_{\varepsilon}\right)  =\mathrm{O}\left(  \varepsilon^{-N_{f}\left(
q,\left\vert \gamma\right\vert \right)  -N_{g}\left(  0,\left\vert
\alpha-\gamma\right\vert \right)  }\right)  \ \mathrm{as}\;\varepsilon
\rightarrow0$. Since $\gamma\leq\alpha$, we get $\left\vert \gamma\right\vert
+\left\vert \alpha-\gamma\right\vert =\left\vert \alpha\right\vert =l$.
According to (\ref{GlimRegD2D2}), there exists $N\in\mathcal{R}^{\prime}$ such
that, for all $k$ and $k^{\prime}\leq k$ in $\mathbb{N}$, $N_{1}\left(
q,k^{\prime}\right)  +N_{2}\left(  0,k-k^{\prime}\right)  \leq$ $N(q,k)$.
Then
\[
\sup_{x\in\Omega}\left\vert \left(  1+\left\vert x\right\vert \right)
^{q}\partial^{\alpha}\left(  f_{\varepsilon}g_{\varepsilon}\right)
(x)\right\vert =\mathrm{O}\left(  \varepsilon^{-N(q,l)}\right)  \ \mathrm{as}%
\;\varepsilon\rightarrow0.
\]
Thus $\mu_{q,l}\left(  f_{\varepsilon}\,g_{\varepsilon}\right)  =\mathrm{O}%
\left(  \varepsilon^{-N(q,l)}\right)  \ \mathrm{as}\;\varepsilon\rightarrow0$
and $\left(  f_{\varepsilon}\,g_{\varepsilon}\right)  _{\varepsilon}%
\in\mathcal{X}_{\mathcal{S}}^{\mathcal{R}^{\prime}}\left(  \Omega\right)  $.

The same kind of estimates shows also that $\mathcal{N}_{\mathcal{S}}\left(
\Omega\right)  $ (or $\mathcal{N}_{\mathcal{S}}^{\mathcal{R}^{\prime}}\left(
\Omega\right)  $) is an ideal of $\mathcal{X}_{\mathcal{S}}^{\mathcal{R}%
^{\prime}}\left(  \Omega\right)  $.\smallskip

\noindent$\left(  b\right)  $~\textit{Presheaf properties}.-~Take $\Omega_{1}$
and $\Omega_{2}$ two open subsets of $\mathbb{R}^{d}$, with $\Omega_{1}%
\subset\Omega_{2}$, $\left(  u_{\varepsilon}\right)  _{\varepsilon}%
\in\mathcal{X}_{\mathcal{S}}^{\mathcal{R}^{\prime}}\left(  \Omega_{2}\right)
$ (\textit{resp}. $\mathcal{N}_{\mathcal{S}}\left(  \Omega_{2}\right)  $). As
$u_{\varepsilon}\in\mathrm{C}^{\infty}\left(  \Omega_{2}\right)  $ for all
$\varepsilon\in\left(  0,1\right]  $, it admits a restriction $u_{\varepsilon
\left\vert \Omega_{1}\right.  }$, and obviously $\left(  u_{\varepsilon
\left\vert \Omega_{1}\right.  }\right)  _{\varepsilon}\in\mathcal{X}%
_{\mathcal{S}}^{\mathcal{R}^{\prime}}\left(  \Omega_{1}\right)  $
(\textit{resp}. $\mathcal{N}_{\mathcal{S}}\left(  \Omega_{2}\right)  $). So,
restrictions are well defined in $\mathcal{X}_{\mathcal{S}}^{\mathcal{R}%
^{\prime}}\left(  \cdot\right)  $ and $\mathcal{N}_{\mathcal{S}}\left(
\cdot\right)  $.
\end{proof}

\begin{definition}
\label{GregDefGSR}The presheaf $\mathcal{G}_{\mathcal{S}}^{\mathcal{R}%
^{\prime}}\left(  \cdot\right)  =\mathcal{X}_{\mathcal{S}}^{\mathcal{R}%
^{\prime}}\left(  \cdot\right)  /\mathcal{N}_{\mathcal{S}}\left(
\cdot\right)  $ is called the presheaf of $\mathcal{R}^{\prime}$\emph{-regular
rapidly decreasing generalized functions}.
\end{definition}

As for the case of $\mathcal{G}^{\mathcal{R}}\left(  \cdot\right)  $, the
presheaf $\mathcal{G}_{\mathcal{S}}^{\mathcal{R}^{\prime}}\left(
\cdot\right)  $ is a presheaf of differential algebras and a sheaf of modules
over the factor ring $\overline{\mathbb{C}}=\mathcal{X}\left(  \mathbb{C}%
\right)  /\mathcal{N}\left(  \mathbb{C}\right)  $.

\begin{example}
\label{GlimSEx1}Taking $\mathcal{R}^{\prime}=\mathbb{R}_{+}^{\mathbb{N}^{2}}$,
we obtain the presheaf of algebras of rapidly decreasing generalized functions
(\cite{Garetto}, \cite{Radyno}, \cite{Scarpa1}).
\end{example}

\begin{notation}
\label{NotaGclassS}In the sequel, we shall note $\mathcal{G}_{\mathcal{S}%
}\left(  \Omega\right)  $ (resp. $\mathcal{X}_{\mathcal{S}}\left(
\Omega\right)  $) instead of $\mathcal{G}_{\mathcal{S}}^{\mathbb{R}%
_{+}^{\mathbb{N}^{2}}}\left(  \Omega\right)  $ (resp. $\mathcal{X}%
_{\mathcal{S}}^{\mathbb{R}_{+}^{\mathbb{N}^{2}}}\left(  \Omega\right)  $). For
all regular subset $\mathcal{R}^{\prime}$ and $\left(  f_{\varepsilon}\right)
_{\varepsilon}\in\mathcal{X}_{\mathcal{S}}^{\mathcal{R}^{\prime}}\left(
\Omega\right)  $, $\left[  \left(  f_{\varepsilon}\right)  _{\varepsilon
}\right]  _{\mathcal{S}}$ $\ $denotes its class in $\mathcal{G}_{\mathcal{S}%
}^{\mathcal{R}^{\prime}}\left(  \Omega\right)  $.
\end{notation}

\begin{example}
\label{GlimSEx2}Taking $\mathcal{R}^{\prime}=\mathcal{B}^{\prime}$, we obtain
the presheaf of $\mathcal{G}_{\mathcal{S}}^{\infty}$ generalized functions or
of \emph{regular rapidly decreasing generalized functions}.
\end{example}

Set
\begin{equation}
\mathcal{N}_{\mathcal{S}_{\ast}}\left(  \Omega\right)  =\left\{  \left(
f_{\varepsilon}\right)  _{\varepsilon}\in\mathrm{C}^{\infty}\left(
\Omega\right)  ^{\left(  0,1\right]  }\,\left\vert \,\forall m\in
\mathbb{N},\;\forall q\in\mathbb{N}\;\;\mu_{q,0}\left(  f_{\varepsilon
}\right)  =\mathrm{O}\left(  \varepsilon^{m}\right)  \;\mathrm{as}%
\;\varepsilon\rightarrow0\right.  \right\}  . \label{GregNSstar}%
\end{equation}
We have the exact analogue of theorems 1.2.25 and 1.2.27 of \cite{GKOS}.

\begin{lemma}
\label{GregLemmaBox}If the open set $\Omega$ is a box, i.e. the product of $d$
open intervals of $\mathbb{R}$ (bounded or not) then $\mathcal{N}%
_{\mathcal{S}}\left(  \Omega\right)  $ is equal to $\mathcal{N}_{\mathcal{S}%
_{\ast}}\left(  \Omega\right)  \cap\mathcal{X}_{\mathcal{S}}\left(
\Omega\right)  $.
\end{lemma}

\begin{proof}
We use similar technics as in the proof of Theorem 1.2.27 in \cite{GKOS}. It
suffices to show the result for $\mathcal{X}_{\mathcal{S}}\left(
\Omega\right)  $ (which is the biggest subspace of $\mathcal{S}\left(
\Omega\right)  ^{\left(  0,1\right]  }$ we may have to consider) and for real
valued functions.\ We prove here that a derivative $\partial_{j}%
=\partial\,/\partial x_{j}$ ($1\leq j\leq d$) satisfies the $0$-order estimate
of the definition of $\mathcal{N}_{\mathcal{S}_{\ast}}$ (\ref{GregNSstar}).
The proof for higher order derivatives, which goes by induction, is left to
the reader. Let $\left(  f_{\varepsilon}\right)  _{\varepsilon}$ be in
$\mathcal{N}_{\mathcal{S}_{\ast}}\left(  \Omega\right)  \cap\mathcal{X}%
_{\mathcal{S}}\left(  \Omega\right)  $, $q$ in $\mathbb{N}$ and $m$ in
$\mathbb{N}$. As $\left(  f_{\varepsilon}\right)  _{\varepsilon}$ is in
$\mathcal{X}_{\mathcal{S}}\left(  \Omega\right)  $, there exists $N$ such that%
\begin{equation}
\forall q^{\prime}\in\left\{  0,...,q\right\}  \ \ \ \ \ \sup_{x\in\Omega
}\left(  1+\left\vert x\right\vert \right)  ^{q\prime}\left\vert \partial
_{j}^{2}f_{\varepsilon}\left(  x\right)  \right\vert =\mathrm{O}\left(
\varepsilon^{-N}\right)  . \label{CSTGEstKern0A}%
\end{equation}
As $\left(  f_{\varepsilon}\right)  _{\varepsilon}$ is in $\mathcal{N}%
_{\mathcal{S}_{\ast}}\left(  \Omega\right)  $, we get
\begin{equation}
\forall q^{\prime}\in\mathbb{N}\ \ \ \ \ \sup_{x\in\Omega}\left(  1+\left\vert
x\right\vert \right)  ^{q\prime}\left\vert f_{\varepsilon}\left(  x\right)
\right\vert =\mathrm{O}\left(  \varepsilon^{N+2m}\right)  .
\label{CSTGEstKern0B}%
\end{equation}
Since the open set $\Omega$ is a box, for $\varepsilon$ sufficiently small
(but independent of $x$), either the segment $\left[  x,x+\varepsilon
^{N+m}e_{j}\right]  $ or $\left[  x,x-\varepsilon^{N+m}e_{j}\right]  $ is
included in $\Omega$. Suppose it is $\left[  x,x+\varepsilon^{N+m}%
e_{j}\right]  $.\ Taylor's theorem gives the existence of $\theta\in\left(
0,1\right)  $ such that
\[
\partial_{j}f_{\varepsilon}\left(  x\right)  =\left(  f_{\varepsilon}\left(
x+\varepsilon^{N+m}e_{j}\right)  -f_{\varepsilon}\left(  x\right)  \right)
\varepsilon^{-N-m}-\left(  1/2\right)  \partial_{j}^{2}f_{\varepsilon}\left(
x_{\theta}\right)  \varepsilon^{N+m},\ \ \ \ \ x_{\theta}=x+\theta
\varepsilon^{N+m}e_{j}.
\]
This gives
\begin{multline*}
\left(  1+\left\vert x\right\vert \right)  ^{q}\left\vert \partial
_{j}f_{\varepsilon}\left(  x\right)  \right\vert \leq\underset{\ast
}{\underbrace{\left(  1+\left\vert x\right\vert \right)  ^{q}\left\vert
f_{\varepsilon}\left(  x+\varepsilon^{N+m}e_{j}\right)  \right\vert
\varepsilon^{-N-m}}}+\underset{\ast\ast}{\underbrace{\left(  1+\left\vert
x\right\vert \right)  ^{q}\left\vert f_{\varepsilon}\left(  x\right)
\right\vert \varepsilon^{-N-m}}}+\\
\underset{\ast\ast\ast}{\underbrace{\varepsilon^{N+m}\left(  1+\left\vert
x\right\vert \right)  ^{q}\left\vert \partial_{j}^{2}f_{\varepsilon}\left(
x_{\theta}\right)  \right\vert }}.
\end{multline*}
From (\ref{CSTGEstKern0B}), we get immediately that $(\ast\ast)$ is of order
$\mathrm{O}\left(  \varepsilon^{m}\right)  $.\ For $\left(  \ast\right)  $, we
have
\[
\left(  1+\left\vert x\right\vert \right)  ^{q}\leq\left(  1+\left\vert
x+\varepsilon^{N+m}e_{j}\right\vert +\varepsilon^{N+m}\right)  ^{q}\leq
\sum_{k=0}^{q}\mathrm{C}_{q}^{k}\left(  1+\left\vert x+\varepsilon^{N+m}%
e_{j}\right\vert \right)  ^{q-k}\varepsilon^{\left(  N+m\right)  k}.
\]
Then
\begin{multline*}
\left(  1+\left\vert x\right\vert \right)  ^{q}\left\vert f_{\varepsilon
}\left(  x+\varepsilon^{N+m}e_{j}\right)  \right\vert \varepsilon^{-N-m}\\
\leq\sum_{k=0}^{q}\mathrm{C}_{q}^{k}\left(  1+\left\vert x+\varepsilon
^{N+m}e_{j}\right\vert \right)  ^{q-k}\left\vert f_{\varepsilon}\left(
x+\varepsilon^{N+m}e_{j}\right)  \right\vert \varepsilon^{\left(  N+m\right)
\left(  k-1\right)  },
\end{multline*}
and (\ref{CSTGEstKern0B}) implies that $\left(  \ast\right)  $ is of order
$\mathrm{O}\left(  \varepsilon^{m}\right)  $. Finally, the same method shows
that $\left(  \ast\ast\ast\right)  $ is also of order $\mathrm{O}\left(
\varepsilon^{m}\right)  $.
\end{proof}

\subsection{Embeddings}

\subsubsection{The natural embeddings of $\mathcal{S}\left(  \mathbb{R}%
^{d}\right)  $ and $\mathcal{O}_{C}^{\prime}\left(  \mathbb{R}^{d}\right)  $
into $\mathcal{G}_{\mathcal{S}}\left(  \mathbb{R}^{d}\right)  $}

The embedding of $\mathcal{S}\left(  \mathbb{R}^{d}\right)  $ into
$\mathcal{G}_{\mathcal{S}}\left(  \mathbb{R}^{d}\right)  $ is done by the
canonical injective map
\[
\sigma_{\mathcal{S}}:\mathcal{S}\left(  \mathbb{R}^{d}\right)  \rightarrow
\mathcal{G}_{\mathcal{S}}\left(  \mathbb{R}^{d}\right)  \;\;\ f\mapsto\left[
\left(  f_{\varepsilon}\right)  _{\varepsilon}\right]  _{\mathcal{S}%
}\ \ \mathrm{with\;}f_{\varepsilon}=f\text{ for all }\varepsilon\in\left(
0,1\right]  .
\]
In fact, the image of $\sigma_{\mathcal{S}}$ is included in $\mathcal{G}%
_{\mathcal{S}}^{\mathcal{R}^{\prime}}\left(  \mathbb{R}^{d}\right)  $ for any
regular subset of $\mathcal{R}^{\prime}\subset\mathbb{R}_{+}^{\mathbb{N}^{2}}$.

For the embedding of $\mathcal{O}_{C}^{\prime}\left(  \mathbb{R}^{d}\right)
$, we consider $\rho\in\mathcal{S}\left(  \mathbb{R}^{d}\right)  $ which
satisfies
\begin{equation}%
{\textstyle\int}
\rho\left(  x\right)  \,\mathrm{d}x=1,\ \ \ \ \ \ \ \ \ \ \
{\textstyle\int}
x^{m}\rho\left(  x\right)  \,\mathrm{d}x=0\;\mathrm{for\;all\;}m\in
\mathbb{N}^{d}\backslash\left\{  0\right\}  \label{CSTGSBinfini}%
\end{equation}
Set
\begin{equation}
\forall\varepsilon\in\left(  0,1\right]  ,\;\;\forall x\in\mathbb{R}%
^{d}\;\;\;\;\rho_{\varepsilon}\left(  x\right)  =\frac{1}{\varepsilon^{d}%
}\,\rho\left(  \frac{x}{\varepsilon}\right)  . \label{CSTGSBinfini2}%
\end{equation}
Note that, contrary to the case of the embedding of $\mathcal{D}^{\prime
}\left(  \Omega\right)  $ in $\mathcal{G}\left(  \Omega\right)  $, we don't
need an additional cutoff.

\begin{theorem}
\label{CSTGThemb}The map
\begin{equation}
\iota_{\mathcal{S}}:\mathcal{O}_{C}^{\prime}\left(  \mathbb{R}^{d}\right)
\rightarrow\mathcal{G}_{\mathcal{S}}\left(  \mathbb{R}^{d}\right)
\;\;\;\;\;u\mapsto\left[  \left(  u\ast\rho_{\varepsilon}\right)
_{\varepsilon}\right]  _{\mathcal{S}}\ \label{CSTGEmbOPCGs}%
\end{equation}
is a linear embedding which commutes with partial derivatives.\newline
\end{theorem}

\begin{proof}
Take $u\in\mathcal{O}_{C}^{\prime}\left(  \mathbb{R}^{d}\right)  $. First, as
$u\in\mathcal{O}_{C}^{\prime}\left(  \mathbb{R}^{d}\right)  $ and
$\rho_{\varepsilon}\in\mathcal{S}\left(  \mathbb{R}^{d}\right)  $,
$u_{\varepsilon}\ast\rho_{\varepsilon}$ is in $\mathcal{S}\left(
\mathbb{R}^{d}\right)  $ for all $\varepsilon\in\left(  0,1\right]  $.\ (This
result is classical: Our proof, with a slight adaptation, shows also this result.)

Consider $q\in\mathbb{N}$. Using the structure of elements of $\mathcal{O}%
_{C}^{\prime}\left(  \mathbb{R}^{d}\right)  $ \cite{Schwartz1}, we can find a
finite family $\left(  f_{j}\right)  _{1\leq j\leq l(q)}$ of continuous
functions such that $\left(  1+\left\vert x\right\vert \right)  ^{q}f_{j}$ is
bounded (for $1\leq j\leq l(q)$), and $\left(  \alpha_{j}\right)  _{1\leq
j\leq l(q)}\in\left(  \mathbb{N}^{d}\right)  ^{l(q)}$ such that $u=\sum
_{j=1}^{l(q)}\partial^{\alpha_{j}}f_{j}$.\ In order to simplify notations, we
shall suppose that this family is reduced to one element $f$, that is
$u=\partial^{\alpha}f$. Take now $\beta\in\mathbb{N}^{d}$. We have
\begin{align*}
\left.  \forall x\in\mathbb{R}^{d}\;\;\;\partial^{\beta}\left(  u\ast
\rho_{\varepsilon}\right)  (x)\right.   &  =\partial^{\beta}\left(
\partial^{\alpha}f\ast\rho_{\varepsilon}\right)  (x)=\left(  f\ast
\partial^{\alpha+\beta}\left(  \rho_{\varepsilon}\right)  \right)  (x),\\
&  =\int f\left(  x-y\right)  \partial^{\alpha+\beta}\left(  \rho
_{\varepsilon}\right)  \left(  y\right)  \,\mathrm{d}y.
\end{align*}
As $\rho_{\varepsilon}\left(  y\right)  =\varepsilon^{-d}\rho\left(
y/\varepsilon\right)  $, we have $\partial^{\alpha+\beta}\left(
\rho_{\varepsilon}\right)  \left(  y\right)  =\varepsilon^{-d-\left\vert
\alpha\right\vert -\left\vert \beta\right\vert }\left(  \partial^{\alpha
+\beta}\rho\right)  \left(  y/\varepsilon\right)  $ and
\[
\forall x\in\mathbb{R}^{d}\;\;\;\partial^{\beta}\left(  u\ast\rho
_{\varepsilon}\right)  (x)=\varepsilon^{-\left\vert \alpha\right\vert
-\left\vert \beta\right\vert }\int f\left(  x-\varepsilon v\right)
\partial^{\alpha+\beta}\rho\left(  v\right)  \,\mathrm{d}v.
\]
On one hand, there exists a constant $C_{1}>0$ such that
\[
\forall\left(  x,v\right)  \in\mathbb{R}^{2d}\;\;\ \left\vert f\left(
x-\varepsilon v\right)  \right\vert \leq C_{1}\left(  1+\left\vert
x-\varepsilon v\right\vert \right)  ^{-q}.
\]
On the other hand, as $\rho$ is rapidly decreasing, there exists $C_{2}>0$
such that
\[
\partial^{\alpha+\beta}\rho\left(  v\right)  \leq C_{2}\left(  1+\left\vert
v\right\vert \right)  ^{-q-d-1}.
\]
These estimates imply the existence of a constant $C_{3}$ such that, for all
$x\in\mathbb{R}^{d}$,
\[
\left.  \left\vert \partial^{\beta}\left(  u\ast\rho_{\varepsilon}\right)
(x)\right\vert \right.  \leq C_{3}\,\varepsilon^{-\left\vert \alpha\right\vert
-\left\vert \beta\right\vert }\int\left(  \left(  1+\left\vert x-\varepsilon
v\right\vert \right)  \left(  1+\left\vert v\right\vert \right)  \right)
^{-q}\left(  1+\left\vert v\right\vert \right)  ^{-d-1}\,\mathrm{d}v.
\]
We have $\left(  1+\left\vert x-\varepsilon v\right\vert \right)  \geq\left(
1+\left\vert \left\vert x\right\vert -\varepsilon\left\vert v\right\vert
\right\vert \right)  $ and a short study of the family of functions
$\phi_{\left\vert x\right\vert ,\varepsilon}:t\mapsto\left(  1+\left\vert
\left\vert x\right\vert -\varepsilon t\right\vert \right)  \left(  1+t\right)
$ for positive $t$ shows that $\phi_{\left\vert x\right\vert ,\varepsilon
}\left(  t\right)  \geq1+\left\vert x\right\vert $.\ Consequently
\begin{align*}
\left.  \forall x\in\mathbb{R}^{d}\;\;\;\left\vert \partial^{\beta}\left(
u\ast\rho_{\varepsilon}\right)  (x)\right\vert \right.   &  \leq
C_{3}\,\varepsilon^{-\left\vert \alpha\right\vert -\left\vert \beta\right\vert
}\left(  1+\left\vert x\right\vert \right)  ^{-q}\int\left(  1+\left\vert
v\right\vert \right)  ^{-d-1}\,\mathrm{d}v\\
&  \leq C_{4}\,\varepsilon^{-\left\vert \alpha\right\vert -\left\vert
\beta\right\vert }\left(  1+\left\vert x\right\vert \right)  ^{-q}\ \text{
(}C_{4}\text{ positive constant).}%
\end{align*}
It follows that $\mu_{q,l}\left(  u\ast\rho_{\varepsilon}\right)
=\mathrm{O}\left(  \varepsilon^{-N(q,l)}\right)  $\ as $\varepsilon
\rightarrow0$ with $N(q,l)=\left\vert \alpha\right\vert +l$. ($\alpha$ may
depends on $q$) This shows that $\left(  u\ast\rho_{\varepsilon}\right)
_{\varepsilon}$ belongs to $\mathcal{X}_{\mathcal{S}}\left(  \mathbb{R}%
^{d}\right)  $.

Finally, it is clear that $\left(  u\ast\rho_{\varepsilon}\right)
_{\varepsilon}\in\mathcal{N}_{\mathcal{S}}\left(  \mathbb{R}^{d}\right)  $
implies that $u_{\varepsilon}\ast\rho_{\varepsilon}\rightarrow0$ in
$\mathcal{S}^{\prime}$, as $\varepsilon\rightarrow0$. As $u_{\varepsilon}%
\ast\rho_{\varepsilon}\rightarrow u$ as $\varepsilon\rightarrow0$, $u$ is
therefore null.
\end{proof}

\begin{theorem}
\label{CSTGThembcom}We have: $\iota_{\mathcal{S}\left\vert \mathcal{S}\left(
\mathbb{R}^{d}\right)  \right.  }=\sigma_{\mathcal{S}}.$
\end{theorem}

\begin{proof}
We shall prove this assertion in the case $d=1$, the general case only differs
by more complicate algebraic expressions.

Let $f$ be in$\mathcal{\ S}\left(  \mathbb{R}\right)  $ and set $\Delta
=\iota_{\mathcal{S}}\left(  f\right)  -\sigma\left(  f\right)  .$ One
representative of $\Delta$ is given by
\[
\Delta_{\varepsilon}:\mathbb{R}\rightarrow\mathcal{F}\left(  \mathrm{C}%
^{\infty}\left(  \mathbb{R}\right)  \right)  \;\;\;\;\;y\mapsto\left(
f\ast\rho_{\varepsilon}\right)  \left(  y\right)  -f(y)=\int f(y-x)\rho
_{\varepsilon}(x)\,\mathrm{d}x-f(y).
\]
Using the fact that $\int\rho_{\varepsilon}\left(  x\right)  \,\mathrm{d}x=1$,
we get
\[
\forall y\in\mathbb{R}\ \ \ \ \Delta_{\varepsilon}(y)=\int\left(
f(y-x)-f(y)\right)  \rho_{\varepsilon}(x)\,\mathrm{d}x.
\]
Let $m$ be an integer. Taylor's formula gives
\[
\forall\left(  x,y\right)  \in\mathbb{R}^{2}\ \ \ \ f(y-x)-f(y)=\sum_{j=1}%
^{m}\frac{\left(  -x\right)  ^{j}}{j!}f^{\left(  i\right)  }\left(  y\right)
+\frac{\left(  -x\right)  ^{m+1}}{m!}\int_{0}^{1}f^{\left(  m+1\right)
}\left(  y-ux\right)  \left(  1-u\right)  ^{m}\,\mathrm{d}u,
\]
and, for all $y\in\mathbb{R}$,
\[
\Delta_{\varepsilon}(y)=\sum_{j=1}^{m}\frac{\left(  -1\right)  ^{j}}%
{j!}f^{\left(  i\right)  }\left(  y\right)  \int x^{i}\rho_{\varepsilon
}(x)\,\mathrm{d}x+\underset{R_{\varepsilon}(m,y)}{\underbrace{\int
\frac{\left(  -x\right)  ^{m+1}}{m!}\left(  \int_{0}^{1}f^{\left(  m+1\right)
}\left(  y-ux\right)  \left(  1-u\right)  ^{m}\mathrm{d}u\right)
\rho_{\varepsilon}(x)\,\mathrm{d}x}}.
\]
According to the choice of mollifiers, we have $\int x^{i}\rho_{\varepsilon
}(x)\,\mathrm{d}x=0$ for $\varepsilon\in\left(  0,1\right]  $ and
$j\in\left\{  0,\ldots,m\right\}  $ and consequently $\Delta_{\varepsilon
}(y)=R_{\varepsilon}(m,y)$.

Setting $v=x/\varepsilon$, we get
\[
\forall y\in\mathbb{R}\ \ \ \ \Delta_{\varepsilon}(y)=\varepsilon^{m+1}%
\int\frac{\left(  -v\right)  ^{m+1}}{m!}\left(  \int_{0}^{1}f^{\left(
m+1\right)  }\left(  y-\varepsilon uv\right)  \left(  1-u\right)
^{m}\,\mathrm{d}u\right)  \,\rho\left(  v\right)  \,\mathrm{d}v.
\]
Consider $q\in\mathbb{N}$.\ As $\rho$ (\textit{resp}. $f$) is in
$\mathcal{S}\left(  \mathbb{R}\right)  $, we get a constant $C_{1}>0$
(\textit{resp}. $C_{2}>0$) such that $\left\vert \rho\left(  t\right)
\right\vert \leq C_{1}\left(  1+\left\vert t\right\vert \right)  ^{-m-3-q}$
(respectively $\left\vert f^{\left(  m+1\right)  }\left(  t\right)
\right\vert \leq C_{2}\left(  1+\left\vert t\right\vert \right)  ^{-q}$) for
all $t\in\mathbb{R}$. Thus, there exists a constant $C_{3}>0$ such that
\[
\forall y\in\mathbb{R}\ \ \ \ \left\vert \Delta_{\varepsilon}(y)\right\vert
\leq C_{3}\frac{\varepsilon^{m+1}}{m!}\int\left\vert v\right\vert
^{m+1}\left(  1+\left\vert v\right\vert \right)  ^{-m-3}\left(  \int_{0}%
^{1}\underline{\left(  1+\left\vert y-\varepsilon uv\right\vert \right)
^{-q}\left(  1+\left\vert v\right\vert \right)  ^{-q}}\,\mathrm{d}u\right)
\mathrm{d}v.
\]
Using the same technic as in the proof of theorem \ref{CSTGThemb}, we obtain
that the underlined term is less than $\left(  1+\left\vert y\right\vert
\right)  ^{-q}$, for all $\varepsilon\in\left(  0,1\right]  $, $u\in\left[
0,1\right]  $ and $v\in\mathbb{R}$. This implies the existence of a constant
$C_{4}>0$ such that
\[
\forall y\in\mathbb{R}\ \ \ \ \left\vert \Delta_{\varepsilon}(y)\right\vert
\leq C_{4}\,\varepsilon^{m+1}\left(  1+\left\vert y\right\vert \right)
^{-q}.
\]
As $\left(  \Delta_{\varepsilon}\right)  _{\varepsilon}\in\mathcal{X}%
_{\mathcal{S}}\left(  \mathbb{R}\right)  $ and $\sup_{y\in\mathbb{R}}\left(
1+\left\vert y\right\vert \right)  ^{q}\left\vert \Delta_{\varepsilon
}(y)\right\vert =\mathrm{O}\left(  \varepsilon^{m+1}\right)  $ for all $m>0$,
we can conclude directly (without estimating the derivatives) that $\left(
\Delta_{\varepsilon}\right)  _{\varepsilon}\in\mathcal{N}_{\mathcal{S}}\left(
\mathbb{R}\right)  $, by using lemma \ref{GregLemmaBox}.\bigskip
\end{proof}

Consider
\begin{equation}
\mathcal{R}_{1}^{\prime}=\left\{  N^{\prime}\in\mathbb{R}_{+}^{\mathbb{N}^{2}%
}\,\left\vert \,\exists N\in\mathcal{R}_{1}\text{ \ }N^{\prime}=1\otimes
N\right.  \right\}  , \label{GregR1PP}%
\end{equation}
where $\mathcal{R}_{1}$ is defined by (\ref{GregDefR1}). (This amounts to:
$N\in\mathcal{R}_{1}^{\prime}$ iff there exists $b\in\mathbb{R}_{+}$ such that
$N(q,l)\leq l+b$.) The set $\mathcal{R}_{1}^{\prime}$ is clearly regular and
we note
\begin{equation}
\mathcal{G}_{\mathcal{S}}^{\left(  1\right)  }\left(  \mathbb{R}^{d}\right)
=\mathcal{X}_{\mathcal{S}}^{\mathcal{R}_{1}^{\prime}}\left(  \mathbb{R}%
^{d}\right)  /\mathcal{N}_{\mathcal{S}}\left(  \mathbb{R}^{d}\right)  .
\label{GregG1SU}%
\end{equation}

\begin{proposition}
\label{GregImbedOPM}The image of $\mathcal{O}_{M}^{\prime}\left(
\mathbb{R}^{d}\right)  $ by $\iota_{\mathcal{S}}$ is included in
$\mathcal{G}_{\mathcal{S}}^{\left(  1\right)  }\left(  \mathbb{R}^{d}\right)
$.
\end{proposition}

\begin{proof}
Let $u$ be in $\mathcal{O}_{M}^{\prime}\left(  \mathbb{R}^{d}\right)  $.
According to the characterization of elements of $\mathcal{O}_{M}^{\prime
}\left(  \mathbb{R}^{d}\right)  $ \cite{Groth}, there exists a finite family
$\left(  f_{j}\right)  _{1\leq j\leq l}$ of rapidly decreasing continuous
functions and $\left(  \alpha_{j}\right)  _{1\leq j\leq l}\in\left(
\mathbb{N}^{d}\right)  ^{l}$ such that $u=\sum_{j=1}^{l}\partial^{\alpha_{j}%
}f_{j}$.\ In order to simplify, we shall suppose that this family is reduced
to one element $f$, that is $u=\partial^{\alpha}f$. For $\beta\in
\mathbb{N}^{d}$, the same estimates as in proof of theorem \ref{CSTGThemb}
lead to the following property
\[
\forall q\in\mathbb{N},\;\;\exists C_{q}>0,\;\;\forall x\in\mathbb{R}%
^{d}\;\;\;\left(  1+\left\vert x\right\vert \right)  ^{q}\left\vert
\partial^{\beta}\left(  u\ast\rho_{\varepsilon}\right)  (x)\right\vert \leq
C_{q}\,\varepsilon^{-\left\vert \alpha\right\vert -\left\vert \beta\right\vert
},
\]
since, in the present case, $f$ is rapidly decreasing. (The only difference is
here that $f$ and $\alpha$ do not depend on the chosen integer $q$.) Then%
\[
\mu_{q,l}\left(  u\ast\rho_{\varepsilon}\right)  \leq C_{q}\,\varepsilon
^{-l-\left\vert \alpha\right\vert }.
\]
Our claim follows, with $N^{\prime}(q,l)=l+\left\vert \alpha\right\vert $,
where $\left\vert \alpha\right\vert $ only depends on $u$.\smallskip\medskip
\end{proof}

We can summarize theorems \ref{CSTGThemb} and \ref{CSTGThembcom} and
proposition \ref{GregImbedOPM} in the following commutative diagram in which
all arrows are embeddings (compare with diagram (\ref{ScomDiagDP})):%

\begin{equation}%
\begin{array}
[c]{ccccc}%
\mathcal{S}\left(  \mathbb{R}^{d}\right)  & \longrightarrow & \mathcal{O}%
_{M}^{\prime}\left(  \mathbb{R}^{d}\right)  & \longrightarrow & \mathcal{O}%
_{C}^{\prime}\left(  \mathbb{R}^{d}\right) \\
& \searrow\sigma_{\mathcal{S}} & \downarrow\iota_{\mathcal{S}} &  &
\downarrow\iota_{\mathcal{S}}\\
&  & \mathcal{G}_{\mathcal{S}}^{(1)}\left(  \mathbb{R}^{d}\right)  &
\longrightarrow & \mathcal{G}_{\mathcal{S}}\left(  \mathbb{R}^{d}\right)  \,.
\end{array}
\label{ScomDiagSP}%
\end{equation}

\subsubsection{Embedding of $\mathcal{S}^{\prime}\left(  \mathbb{R}%
^{d}\right)  $ into $\mathcal{G}_{\mathcal{S}}\left(  \mathbb{R}^{d}\right)
$}

In order to embed $\mathcal{S}^{\prime}\left(  \mathbb{R}^{d}\right)  $ into
an algebra playing the role of $\mathcal{G}\left(  \mathbb{R}^{d}\right)  $
for $\mathcal{D}^{\prime}\left(  \mathbb{R}^{d}\right)  $, a space
$\mathcal{G}_{\tau}\left(  \mathbb{R}^{d}\right)  $ of tempered generalized
functions is often introduced (see \cite{Col1}, \cite{GKOS}). This space
$\mathcal{G}_{\tau}\left(  \mathbb{R}^{d}\right)  $ does not fit in the
general scheme of construction of Colombeau type algebras, since the growth
estimates for $\mathcal{G}_{\tau}\left(  \mathbb{R}^{d}\right)  $ are not
based on the natural topology of the space $\mathcal{O}_{M}\left(
\mathbb{R}^{d}\right)  $, which replaces $\mathrm{C}^{\infty}\left(
\mathbb{R}^{d}\right)  $ in this case. Although it is possible to construct a
space $\mathcal{G}_{\tau}\left(  \mathbb{R}^{d}\right)  $ based on the
topology of $\mathcal{O}_{M}\left(  \mathbb{R}^{d}\right)  $, we don't need it
in the sequel and we only show that $\mathcal{S}^{\prime}\left(
\mathbb{R}^{d}\right)  $ can be embedded into $\mathcal{G}_{\mathcal{S}%
}\left(  \mathbb{R}^{d}\right)  $ by means of a cutoff\ of the embedding
$\iota_{\mathcal{S}}:\mathcal{O}_{C}^{\prime}\left(  \mathbb{R}^{d}\right)
\rightarrow\mathcal{G}_{\mathcal{S}}\left(  \mathbb{R}^{d}\right)  $.

\begin{proposition}
\label{RDGFembSPr}With the notations (\ref{CSTGSBinfini}) and
(\ref{CSTGSBinfini2}), the map
\[
\iota_{\mathcal{S}^{\prime}}:\mathcal{S}^{\prime}\left(  \mathbb{R}%
^{d}\right)  \rightarrow\mathcal{G}_{\mathcal{S}}\left(  \mathbb{R}%
^{d}\right)  ,\;u\mapsto\left[  \left(  \left(  u\ast\rho_{\varepsilon
}\right)  \widehat{\rho}_{\varepsilon}\right)  _{\varepsilon}\right]
_{\mathcal{S}}%
\]
is a linear embedding.
\end{proposition}

\begin{proof}
For $u\in\mathcal{S}^{\prime}\left(  \mathbb{R}^{d}\right)  $, the net
$\left(  \left(  u\ast\rho_{\varepsilon}\right)  \right)  _{\varepsilon}$
belongs to the space%
\[
\mathcal{X}_{\tau}\left(  \mathbb{R}^{d}\right)  =\left\{  \left(
f_{\varepsilon}\right)  _{\varepsilon}\in\mathrm{C}^{\infty}\left(
\Omega\right)  ^{\left(  0,1\right]  }\,\left\vert \,\forall\alpha
\in\mathbb{N}^{d},\;\exists q\in\mathbb{N},\;\exists N\in\mathbb{N}%
\;\;\mu_{-q,\alpha}\left(  f_{\varepsilon}\right)  =\mathrm{O}\left(
\varepsilon^{-N}\right)  \;\mathrm{as}\;\varepsilon\rightarrow0\right.
\right\}  .
\]
(See, for example, theorem 1.2.27 in \cite{GKOS}: The proof is similar to the
one of theorem \ref{CSTGThemb}.) A straightforward calculus shows that
$\mathcal{X}_{S}$ is an ideal of $\mathcal{X}_{\tau}$.\ It follows that
$\left(  \left(  u\ast\rho_{\varepsilon}\right)  \widehat{\rho}_{\varepsilon
}\right)  _{\varepsilon}\in\mathcal{X}_{\mathcal{S}}\left(  \mathbb{R}%
^{d}\right)  $ since $\left(  \widehat{\rho}_{\varepsilon}\right)
_{\varepsilon}\in\mathcal{X}_{\mathcal{S}}\left(  \mathbb{R}^{d}\right)  $.
Then, the map $\iota_{\mathcal{S}^{\prime}}$ is well defined.

Note that, for $u\in\mathcal{S}^{\prime}\left(  \mathbb{R}^{d}\right)  $, we
have $u\ast\rho_{\varepsilon}\overset{\mathcal{S}^{\prime}}{\longrightarrow}u$
as $\varepsilon\rightarrow0$ and, therefore, $u\ast\rho_{\varepsilon}%
\overset{\mathcal{D}^{\prime}}{\longrightarrow}u$ as $\varepsilon\rightarrow
0$. As $\widehat{\rho}_{\varepsilon}=1$ on a compact set $K_{\varepsilon}$
such that $K_{\varepsilon}\rightarrow\mathbb{R}^{d}$ as $\varepsilon
\rightarrow0$, we get that%
\begin{equation}
\left(  u\ast\rho_{\varepsilon}\right)  \widehat{\rho}_{\varepsilon}%
\overset{\mathcal{D}^{\prime}}{\longrightarrow}u\ \mathrm{as\ }\varepsilon
\rightarrow0. \label{RDGFembSprDpr}%
\end{equation}
Finally, take $u\in\mathcal{S}^{\prime}\left(  \mathbb{R}^{d}\right)  $ with
$\iota_{\mathcal{S}^{\prime}}\left(  u\right)  =0$, that is $\left(  \left(
u\ast\rho_{\varepsilon}\right)  \widehat{\rho}_{\varepsilon}\right)
_{\varepsilon}\in\mathcal{N}_{\mathcal{S}}\left(  \mathbb{R}^{d}\right)
$.\ We have $\left(  u\ast\rho_{\varepsilon}\right)  \widehat{\rho
}_{\varepsilon}\overset{\mathcal{S}^{\prime}}{\longrightarrow}0$ and,
consequently, $\left(  u\ast\rho_{\varepsilon}\right)  \widehat{\rho
}_{\varepsilon}\overset{\mathcal{D}^{\prime}}{\longrightarrow}0$.\ Thus, $u=0$
according to (\ref{RDGFembSprDpr}).
\end{proof}

\section{Fourier transform and exchange theorem\label{Glimsec4}}

\subsection{Fourier transform in $\mathcal{G}_{\mathcal{S}}\left(
\mathbb{R}^{d}\right)  $}

The Fourier transform $\mathcal{F}$ is a continuous linear map from
$\mathcal{S}\left(  \mathbb{R}^{d}\right)  $ to $\mathcal{S}\left(
\mathbb{R}^{d}\right)  $.\ According to proposition 3.2.\ of \cite{ADDS},
$\mathcal{F}$ has a canonical extension $\mathcal{F}_{\mathcal{S}}$ from
$\mathcal{G}_{\mathcal{S}}$ to $\mathcal{G}_{\mathcal{S}}$ defined by
\begin{equation}
\mathcal{G}_{\mathcal{S}}\left(  \mathbb{R}^{d}\right)  \rightarrow
\mathcal{G}_{\mathcal{S}}\left(  \mathbb{R}^{d}\right)  \;\;\;u\mapsto
\widehat{u}=\left[  \left(  x\mapsto\int\mathrm{e}^{-ix\xi}u_{\varepsilon
}\left(  \xi\right)  \,\mathrm{d}\xi\right)  _{\varepsilon}\right]
_{\mathcal{S}}\,, \label{GregFTForm}%
\end{equation}
where $\left(  u_{\varepsilon}\right)  _{\varepsilon}\in\mathcal{X}%
_{\mathcal{S}}\left(  \mathbb{R}^{d}\right)  $ is any representative of
$u$.\smallskip

The proof of this result uses mainly the continuity of $\mathcal{F}$.\ More
precisely, any linear continuous map is continuously moderate in the sense of
\cite{ADDS} and, therefore, admits such a canonical extension. (The proof of
lemma \ref{GregLemFond} below, with a slight adaptation, shows directly this
result for the case of $\mathcal{F}$.)

\begin{definition}
\label{GregDefTFGS}The map $\mathcal{F}_{\mathcal{S}}$ defined by
(\ref{GregFTForm}) is called the Fourier transform in $\mathcal{G}%
_{\mathcal{S}}$.
\end{definition}

In the same way, we can define $\mathcal{F}_{\mathcal{S}}^{-1}$ by
\begin{equation}
\mathcal{G}_{\mathcal{S}}\left(  \mathbb{R}^{d}\right)  \rightarrow
\mathcal{G}_{\mathcal{S}}\left(  \mathbb{R}^{d}\right)  \;\;\;u\mapsto\left[
\left(  x\mapsto\left(  2\pi\right)  ^{-d}\int\mathrm{e}^{ix\xi}%
u_{\varepsilon}\left(  \xi\right)  \,\mathrm{d}\xi\right)  _{\varepsilon
}\right]  _{\mathcal{S}}\,, \label{CSTGDefInvF}%
\end{equation}
where $\left(  u_{\varepsilon}\right)  _{\varepsilon}\in\mathcal{X}%
_{\mathcal{S}}\left(  \mathbb{R}^{d}\right)  $ is any representative of $u$.

\begin{theorem}
\label{GregThFourier}$\mathcal{F}_{\mathcal{S}}:\mathcal{G}_{\mathcal{S}%
}\left(  \mathbb{R}^{d}\right)  \rightarrow\mathcal{G}_{_{\mathcal{S}}}\left(
\mathbb{R}^{d}\right)  $ is a one to one linear map, whose inverse is
$\mathcal{F}_{\mathcal{S}}^{-1}$.
\end{theorem}

\begin{proof}
Let $u$ be in $\mathcal{G}_{\mathcal{S}}\left(  \mathbb{R}^{d}\right)  $ and
$\left(  u_{\varepsilon}\right)  _{\varepsilon}\in\mathcal{X}_{\mathcal{S}%
}\left(  \mathbb{R}^{d}\right)  $ be one of its representative. As
representative of $\mathcal{F}\left(  \mathcal{F}^{-1}\left(  u\right)
\right)  $, we can choose $\left(  \widetilde{u}_{\varepsilon}\right)
_{\varepsilon}$ defined by
\[
\forall\varepsilon\in\left(  0,1\right]  ,\;\;\forall x\in\mathbb{R}%
^{d}\;\;\;\widetilde{u}_{\varepsilon}(x)=\left(  2\pi\right)  ^{-d}%
\int\mathrm{e}^{ix\xi}\widehat{u}_{\varepsilon}\left(  \xi\right)
\,\mathrm{d}\xi.
\]
Since the Fourier transform is an isomorphism in $\mathcal{S}\left(
\mathbb{R}^{d}\right)  $, we get $\widetilde{u}_{\varepsilon}=u_{\varepsilon}%
$, for all $\varepsilon\in\left(  0,1\right]  $, and $\mathcal{F}%
_{\mathcal{S}}\left(  \mathcal{F}_{\mathcal{S}}^{-1}\left(  u\right)  \right)
=\left[  \left(  u_{\varepsilon}\left(  x\right)  \right)  _{\varepsilon
}\right]  _{\mathcal{S}}=u$.
\end{proof}

\subsection{Regular sub presheaves of $\mathcal{G}_{\mathcal{S}}\left(
\cdot\right)  $}

We introduce here some regular sub presheaves of $\mathcal{G}_{\mathcal{S}%
}\left(  \cdot\right)  $ needed for our further microlocal analysis.

Let $\mathcal{R}$ be a regular subset of $\mathbb{R}_{+}^{\mathbb{N}}$ and set%
\[
\mathcal{R}_{u}=\left\{  N^{\prime}\in\mathbb{R}_{+}^{\mathbb{N}^{2}%
}\,\left\vert \,\exists N\in\mathcal{R}\ \ \ N^{\prime}=1\otimes N\right.
\right\}  ;\ \ \ \mathcal{R}_{\partial}=\left\{  N^{\prime}\in\mathbb{R}%
_{+}^{\mathbb{N}^{2}}\,\left\vert \,\exists N\in\mathcal{R}\ \ \ N=N\otimes
1\right.  \right\}  .
\]
In other words, $N^{\prime}\in\mathcal{R}_{u}$ (\textit{resp}. $\mathcal{R}%
_{\partial}$) iff there exists $N\in\mathcal{R}$ such that $N^{\prime
}(q,l)=N(l)$ (\textit{resp}. $N^{\prime}(q,l)=N(q)$) or, equivalently, iff $N$
only depends (in a $\mathcal{R}$-regular way) of $l$ (\textit{resp}. $q$).

\begin{notation}
We shall write, with a slight abuse, $\mathcal{R}_{u}=\left\{  1\right\}
\otimes\mathcal{R}$,\ $\mathcal{R}_{\partial}=\mathcal{R}\otimes\left\{
1\right\}  .$
\end{notation}

Obviously, $\mathcal{R}_{u}$ (\textit{resp}. $\mathcal{R}_{\partial}$) is a
regular subset of $\mathbb{R}_{+}^{\mathbb{N}^{2}}$.

\begin{example}
\label{GregexampDU}Take $\mathcal{R=}\mathbb{R}_{+}^{\mathbb{N}}$. We set:
$\mathcal{G}_{\mathcal{S}}^{u}\left(  \cdot\right)  =\mathcal{G}_{\mathcal{S}%
}^{\mathcal{R}_{u}}\left(  \cdot\right)  $ (resp. $\mathcal{G}_{\mathcal{S}%
}^{\partial}\left(  \cdot\right)  =\mathcal{G}_{\mathcal{S}}^{\mathcal{R}%
_{\partial}}\left(  \cdot\right)  $). In this case, we have $\mathcal{R}%
_{u}=\left\{  1\right\}  \otimes\mathbb{R}_{+}^{\mathbb{N}}$ (resp.
$\mathcal{R}_{\partial}=\mathbb{R}_{+}^{\mathbb{N}}\otimes\left\{  1\right\}
$).
\end{example}

The elements of $\mathcal{G}_{\mathcal{S}}^{u}\left(  \Omega\right)  $
($\Omega$ open subset of $\mathbb{R}^{d}$) have uniform growth bounds with
respect to the regularization parameter $\varepsilon$ for all factors $\left(
1+\left\vert x\right\vert \right)  ^{q}$. For $\mathcal{G}_{\mathcal{S}%
}^{\partial}\left(  \Omega\right)  $, those bounds are uniform for all
derivatives.\ For $\mathcal{G}_{\mathcal{S}}^{\infty}\left(  \Omega\right)  $,
introduced in example \ref{GlimSEx2}, the uniformity is global, in some sense
stronger than the $\mathcal{G}^{\infty}$-regularity considered for the algebra
$\mathcal{G}$. (In this last case, the uniformity is not required with respect
to the compact sets.)

We have the obvious inclusions (for $\mathcal{R}$ regular subset of
$\mathbb{R}_{+}^{\mathbb{N}}$):%

\begin{equation}%
\begin{array}
[c]{ccccccc}
&  & \mathcal{G}_{\mathcal{S}}^{\mathcal{R}_{\partial}}\left(  \Omega\right)
& \longrightarrow & \mathcal{G}_{\mathcal{S}}^{\partial}\left(  \Omega\right)
&  & \\
\mathcal{G}_{\mathcal{S}}^{\infty}\left(  \Omega\right)  &
\begin{array}
[c]{c}%
\nearrow\\
\searrow
\end{array}
&  &  &  &
\begin{array}
[c]{c}%
\searrow\\
\nearrow
\end{array}
& \mathcal{G}_{\mathcal{S}}\left(  \Omega\right) \\
&  & \mathcal{G}_{\mathcal{S}}^{\mathcal{R}_{u}}\left(  \Omega\right)  &
\longrightarrow & \mathcal{G}_{\mathcal{S}}^{u}\left(  \Omega\right)  &  &
\end{array}
. \label{GregDiagramDU1}%
\end{equation}

\begin{example}
\label{GregGR1PS}The algebra $\mathcal{G}_{\mathcal{S}}^{\left(  1\right)
}\left(  \mathbb{R}^{d}\right)  =\mathcal{X}_{\mathcal{S}}^{\mathcal{R}%
_{1}^{\prime}}\left(  \mathbb{R}^{d}\right)  /\mathcal{N}_{\mathcal{S}}\left(
\mathbb{R}^{d}\right)  $ introduced in (\ref{GregG1SU}) for the embedding of
$\mathcal{O}_{M}^{\prime}\left(  \mathbb{R}^{d}\right)  $ into $\mathcal{G}%
_{\mathcal{S}}\left(  \mathbb{R}^{d}\right)  $ (proposition \ref{GregImbedOPM}%
) can be written as $\mathcal{G}_{\mathcal{S}}^{\left(  \mathcal{R}%
_{1}\right)  _{u}}$, with
\[
\mathcal{R}_{1}=\left\{  N\in\mathbb{R}_{+}^{\mathbb{N}}\,\left\vert \,\exists
b\in\mathbb{R}_{+},\ \ \forall l\in\mathbb{R}\ \ \ N(l)\leq l+b\right.
\right\}  .
\]

\end{example}

As a first illustration of the properties of these spaces, we can show the
existence of a canonical embedding of algebras of compactly supported
generalized functions into particular spaces of rapidly decreasing generalized functions.

\begin{lemma}
\label{GregEmbCRSRu}Let $\mathcal{R}$ be a regular subset of $\mathbb{R}%
_{+}^{\mathbb{N}}$ and $u$ be in $\mathcal{G}_{C}^{\mathcal{R}}\left(
\Omega\right)  $ ($\Omega$ open subset of $\mathbb{R}^{d}$), with $\left(
u_{\varepsilon}\right)  _{\varepsilon}$ a representative of $u$. Let $\kappa$
be in $\mathcal{D}\left(  \Omega\right)  $, with $0\leq\kappa\leq1$ and
$\kappa\equiv1$ on a neighborhood of $\operatorname*{supp}u.$ Then $\left(
\kappa u_{\varepsilon}\right)  _{\varepsilon}$ belongs to $\mathcal{X}%
_{\mathcal{S}}^{\mathcal{R}_{u}}\left(  \mathbb{R}^{d}\right)  $ and $\left[
\left(  \kappa u_{\varepsilon}\right)  _{\varepsilon}\right]  _{\mathcal{S}}$
only depends on $u$, but not on $\left(  u_{\varepsilon}\right)
_{\varepsilon}$ and $\kappa$.
\end{lemma}

\begin{proof}
We first show that $\left(  \kappa u_{\varepsilon}\right)  _{\varepsilon}$ is
in $\mathcal{X}_{\mathcal{S}}^{\mathcal{R}_{u}}\left(  \mathbb{R}^{d}\right)
$ and then the independence with respect to the representation.

\noindent$\left(  a\right)  ~$There exists a compact set $K\subset\Omega$ such
that, for all $\varepsilon\in\left(  0,1\right]  $, $\operatorname*{supp}%
\kappa u_{\varepsilon}\subset K$. It follows that $\kappa u_{\varepsilon}$ is
compactly supported and therefore rapidly decreasing. Moreover%
\[
\forall\left(  q,l\right)  \in\mathbb{N}^{2},\ \forall\varepsilon\in\left(
0,1\right]  \ \ \ \mu_{q,l}\left(  \kappa u_{\varepsilon}\right)  \leq
\sup\nolimits_{x\in K}\left(  1+\left\vert x\right\vert \right)  ^{q}%
p_{K,l}\left(  \kappa u_{\varepsilon}\right)  \leq C_{K,q}\,p_{K,l}\left(
\kappa u_{\varepsilon}\right)  ,\ \ \ C_{K,q}>0\text{.}%
\]
Thus, $\left(  \kappa u_{\varepsilon}\right)  _{\varepsilon}$ belongs to
$\mathcal{X}_{\mathcal{S}}^{\mathcal{R}_{u}}\left(  \mathbb{R}^{d}\right)  $.
Indeed, by using Leibniz' formula for estimating $p_{K,l}\left(  \kappa
u_{\varepsilon}\right)  $, we can find a constant $C_{K,q,\kappa}>0$ such that%
\begin{equation}
\forall\left(  q,l\right)  \in\mathbb{N}^{2},\ \ \forall\varepsilon\in\left(
0,1\right]  \ \ \ \ \mu_{q,l}\left(  \kappa u_{\varepsilon}\right)  \leq
C_{K,q,\kappa}p_{K,l}\left(  u_{\varepsilon}\right)  . \label{Glimeq6}%
\end{equation}

\noindent$\left(  b\right)  ~$Let $\left(  \widetilde{u}_{\varepsilon}\right)
_{\varepsilon}$ be another representative of $u$ and $\widetilde{\kappa}$ be
in $\mathcal{D}\left(  \Omega\right)  $, with $0\leq\widetilde{\kappa}\leq1$
and $\widetilde{\kappa}=1$ on a neighborhood of $\operatorname*{supp}%
\widetilde{u}$. Let $L$ be a compact set such that $\operatorname*{supp}\kappa
u_{\varepsilon}\cup\operatorname*{supp}\widetilde{\kappa}\widetilde
{u}_{\varepsilon}\subset L\subset\Omega$.\ According to the previous estimate,
we have
\begin{align*}
\left.  \ \forall\left(  q,l\right)  \in\mathbb{N}^{2},\ \ \forall
\varepsilon\in\left(  0,1\right]  \ \ \ \ \right.  \mu_{q,l}\left(  \kappa
u_{\varepsilon}-\widetilde{\kappa}\widetilde{u}_{\varepsilon}\right)   &
\leq\mu_{q,l}\left(  \left(  \kappa-\widetilde{\kappa}\right)  u_{\varepsilon
}\right)  +\mu_{q,l}\left(  \widetilde{\kappa}\left(  u_{\varepsilon
}-\widetilde{u}_{\varepsilon}\right)  \right) \\
&  \leq C_{L,q}\,p_{L,l}\left(  \left(  \kappa-\widetilde{\kappa}\right)
u_{\varepsilon}\right)  +C_{L,q}\,p_{L,l}\left(  \widetilde{\kappa}\left(
u_{\varepsilon}-\widetilde{u}_{\varepsilon}\right)  \right)  .
\end{align*}
As $\kappa=\widetilde{\kappa}$ on a closed neighborhood $V$ of
$\operatorname*{supp}u$, it follows that $p_{V,l}\left(  \left(
\kappa-\widetilde{\kappa}\right)  u_{\varepsilon}\right)  =0$. Moreover, for
all $m\in\mathbb{N}$, $p_{L\backslash V,l}\left(  \left(  \kappa
-\widetilde{\kappa}\right)  u_{\varepsilon}\right)  =\mathrm{O}\left(
\varepsilon^{m}\right)  $ as $\varepsilon\rightarrow0$, since $\left(
L\backslash V\right)  \cap\operatorname*{supp}u=\varnothing$. Then
$p_{L,l}\left(  \left(  \kappa-\widetilde{\kappa}\right)  u_{\varepsilon
}\right)  =\mathrm{O}\left(  \varepsilon^{m}\right)  $ as $\varepsilon
\rightarrow0$. As $\left[  \left(  u_{\varepsilon}\right)  _{\varepsilon
}\right]  =\left[  \left(  \widetilde{u}_{\varepsilon}\right)  _{\varepsilon
}\right]  $, we have $p_{L,l}\left(  \widetilde{\kappa}\left(  u_{\varepsilon
}-\widetilde{u}_{\varepsilon}\right)  \right)  =\mathrm{O}\left(
\varepsilon^{m}\right)  $ as $\varepsilon\rightarrow0$.\ Then $\mu
_{q,l}\left(  \kappa u_{\varepsilon}-\widetilde{\kappa}\widetilde
{u}_{\varepsilon}\right)  =\mathrm{O}\left(  \varepsilon^{m}\right)  $ and
$\left[  \left(  \kappa u_{\varepsilon}\right)  _{\varepsilon}\right]
_{\mathcal{S}}=\left[  \left(  \widetilde{\kappa}\widetilde{u}\right)
_{\varepsilon}\right]  _{\mathcal{S}}$. \medskip\smallskip
\end{proof}

From lemma \ref{GregEmbCRSRu}, we deduce easily the following proposition:

\begin{proposition}
\label{GregEmbGCRGSRU}With the notations of lemma \ref{GregEmbCRSRu}, the map
\[
\iota_{C,S}:\mathcal{G}_{C}^{\mathcal{R}}\left(  \Omega\right)  \rightarrow
\mathcal{G}_{\mathcal{S}}^{\mathcal{R}_{u}}\left(  \mathbb{R}^{d}\right)
\ \ \ \ \ \ u\mapsto\left[  \left(  \kappa u_{\varepsilon}\right)
_{\varepsilon}\right]  _{\mathcal{S}}\
\]
is a linear embedding.
\end{proposition}

From the embedding $\iota_{C,S}$, one can then verify that the Fourier
transform of a compactly supported generalized functions $u\in\mathcal{G}%
^{\mathcal{R}}\left(  \Omega\right)  $, which can be straightforwardly
considered as an element of $\mathcal{G}_{C}\left(  \mathbb{R}^{d}\right)  $,
is defined by one of the following equalities
\[
\mathcal{F}\left(  u\right)  =\mathcal{F}\left(  \iota_{C,S}(u)\right)
=\left[  \left(  x\mapsto\left(  2\pi\right)  ^{-d}\int_{W}\mathrm{e}^{ix\xi
}u_{\varepsilon}\left(  \xi\right)  \,\mathrm{d}\xi\right)  _{\varepsilon
}\right]  _{\mathcal{S}},
\]
where $\left(  u_{\varepsilon}\right)  _{\varepsilon}\in\mathcal{X}%
^{\mathcal{R}}\left(  \mathbb{R}^{d}\right)  $ is any representative of $u$
and $W$ any relatively compact neighborhood of $\operatorname*{supp}u$.

\subsection{Exchange and regularity theorems}

\begin{theorem}
\label{GregExchangethm}\textbf{(Exchange theorem)} For any regular subset
$\mathcal{R}$ of $\mathbb{R}_{+}^{\mathbb{N}}$, we have:
\begin{equation}
\mathcal{F}\left(  \mathcal{G}_{\mathcal{S}}^{\mathcal{R}_{u}}\left(
\mathbb{R}^{d}\right)  \right)  =\mathcal{G}_{\mathcal{S}}^{\mathcal{R}%
_{\partial}}\left(  \mathbb{R}^{d}\right)
\ \ \ \ \ \ \ \ \ \ \ \ \ \ \ \mathcal{F}\left(  \mathcal{G}_{\mathcal{S}%
}^{\mathcal{R}_{\partial}}\left(  \mathbb{R}^{d}\right)  \right)
=\mathcal{G}_{\mathcal{S}}^{\mathcal{R}_{u}}\left(  \mathbb{R}^{d}\right)  .
\label{GregExcProp}%
\end{equation}

\end{theorem}

\begin{proof}
The proof of this theorem is based on the following classical:

\begin{lemma}
\label{GregLemFond}For all $u\in\mathcal{S}\left(  \mathbb{R}^{d}\right)  $
and $\left(  q,l\right)  \in\mathbb{N}^{2}$, there exists a constant
$C_{q,l}>0$ such that
\begin{equation}
\mu_{q,l}\left(  \widehat{u}\right)  \leq C_{q,l}\,\mu_{l+d+1,q}\left(
u\right)  . \label{GReglimEq9}%
\end{equation}

\end{lemma}

Let us prove this result first. Take $u\in\mathcal{S}\left(  \mathbb{R}%
^{d}\right)  $, $\left(  q,l\right)  \in\mathbb{N}^{2}$ and $\left(
\alpha,\beta\right)  \in\left(  \mathbb{N}^{d}\right)  ^{2}$ with $\left\vert
\alpha\right\vert =l$ and $\left\vert \beta\right\vert \leq q$.\ We have
\[
\forall\xi\in\mathbb{R}^{d}\;\;\;\partial^{\alpha}\widehat{u}\left(
\xi\right)  =\int\mathrm{e}^{-ix\xi}\left(  -ix\right)  ^{\alpha}u\left(
x\right)  \,\mathrm{d}x.
\]
By integration by parts, we obtain
\begin{equation}
\forall\xi\in\mathbb{R}^{d}\;\;\;\xi^{\beta}\partial^{\alpha}\widehat
{u}\left(  \xi\right)  =\int\xi^{\beta}\mathrm{e}^{-ix\xi}\left(  -ix\right)
^{\alpha}u\left(  x\right)  \,\mathrm{d}x=\left(  -i\right)  ^{\left\vert
\beta\right\vert }\int\mathrm{e}^{ix\xi}\partial^{\beta}\left[  \left(
-ix\right)  ^{\alpha}u\left(  x\right)  \right]  \,\mathrm{d}x.
\label{GregGLMTF0}%
\end{equation}
By the Leibniz formula, we have%
\begin{equation}
\forall x\in\mathbb{R}^{d}\;\;\;\partial^{\beta}\left[  \left(  -ix\right)
^{\alpha}u\left(  x\right)  \right]  =\sum_{\gamma\leq\beta}C_{_{\beta}%
}^{\gamma}\partial^{\beta-\gamma}\left[  \left(  -ix\right)  ^{\alpha}\right]
\partial^{\gamma}u\left(  x\right)  , \label{GregGprLmTF1}%
\end{equation}
where$\;C_{_{\beta}}^{\gamma}$ are the $d$ dimensional binomial coefficients.

There exists a constant $C_{\alpha,\beta}^{\prime}>0$ such that, for all
$\gamma\in\mathbb{N}^{d}$ with $\gamma\leq\beta$,
\begin{multline*}
\forall x\in\mathbb{R}^{d}\;\;\;\left(  1+\left\vert x\right\vert \right)
^{d+1}\left\vert \partial^{\beta-\gamma}\left[  \left(  -ix\right)  ^{\alpha
}\right]  \partial^{\gamma}u\left(  x\right)  \right\vert \leq C_{\alpha
,\beta}^{\prime}\left(  1+\left\vert x\right\vert \right)  ^{\left\vert
\alpha\right\vert +d+1}\left\vert \partial^{\gamma}u\left(  x\right)
\right\vert \\
\leq C_{\alpha,\beta}^{\prime}\,\mu_{\left\vert \alpha\right\vert +d+1,\gamma
}\left(  u\right)  \leq C_{\alpha,\beta}^{\prime}\,\sup_{\gamma\leq\beta}%
\mu_{\left\vert \alpha\right\vert +d+1,\gamma}\left(  u\right)  .
\end{multline*}
Summing up these results in (\ref{GregGprLmTF1}), we find a constant
$C_{\alpha,\beta}^{\prime\prime}>0$ such that
\[
\forall x\in\mathbb{R}^{d}\;\;\;\left(  1+\left\vert x\right\vert \right)
^{d+1}\left\vert \partial^{\beta}\left[  \left(  -ix\right)  ^{\alpha}u\left(
x\right)  \right]  \right\vert \leq C_{\alpha,\beta}^{\prime\prime}%
\,\sup_{\gamma\leq\beta}\mu_{\left\vert \alpha\right\vert +d+1,\gamma}\left(
u\right)  .
\]
Going back to relation (\ref{GregGLMTF0}), we have
\[
\forall\xi\in\mathbb{R}^{d}\;\;\;\left\vert \xi^{\beta}\partial^{\alpha
}\widehat{u}\left(  \xi\right)  \right\vert \leq C_{\alpha,\beta}%
^{\prime\prime}\,\sup_{\gamma\leq\beta}\mu_{\left\vert \alpha\right\vert
+d+1,\gamma}\left(  u\right)  \int\left(  1+\left\vert x\right\vert \right)
^{-d-1}\,\mathrm{d}x.
\]
Finally, we get the existence of a constant $C_{\alpha,\beta}>0$
\[
\forall\xi\in\mathbb{R}^{d}\;\;\;\left\vert \xi^{\beta}\partial^{\alpha
}\widehat{u}\left(  \xi\right)  \right\vert \leq C_{\alpha,\beta}%
\,\mu_{\left\vert \alpha\right\vert +d+1,\left\vert \beta\right\vert }\left(
u\right)  \leq C_{q,l}\,\mu_{l+d+1,q}\left(  u\right)  ,
\]
where $C_{q,l}$ is a constant greater than all $C_{\alpha,\beta}$ for
$\left\vert \alpha\right\vert =l$ and $\left\vert \beta\right\vert \leq q$. In
the classical manner, we can deduce inequality (\ref{GReglimEq9}) from this
last estimate.\medskip

We return to the proof of theorem \ref{GregExchangethm}. Let $\mathcal{R}$ be
a regular subset of $\mathbb{R}_{+}^{\mathbb{N}}$\textit{.}

\noindent$\left(  a\right)  $~Take $u\in\mathcal{G}_{\mathcal{S}}%
^{\mathcal{R}_{u}}\left(  \mathbb{R}^{d}\right)  $ and $\left(  u_{\varepsilon
}\right)  _{\varepsilon}\in\mathcal{X}_{\mathcal{S}}^{\mathcal{R}_{u}}\left(
\mathbb{R}^{d}\right)  $ a representative of $u$. There exists a sequence
$N\in\mathcal{R}$ such that $\mu_{r,q}\left(  u_{\varepsilon}\right)
=\mathrm{O}\left(  \varepsilon^{-N(q)}\right)  \;\mathrm{as}\;\varepsilon
\rightarrow0$, for all $r\in\mathbb{N}$.\ Lemma \ref{GregLemFond} implies that
$\mu_{q,l}\left(  \widehat{u}_{\varepsilon}\right)  =\mathrm{O}\left(
\varepsilon^{-N(q)}\right)  \;\mathrm{as}\;\varepsilon\rightarrow0$, for all
$l\in\mathbb{N}$.\ Thus, $\mathcal{F}\left(  u\right)  \in\mathcal{G}%
_{\mathcal{S}}^{\mathcal{R}_{\partial}}\left(  \mathbb{R}^{d}\right)  $.

\noindent$\left(  b\right)  $~Conversely, take $u\in\mathcal{G}_{\mathcal{S}%
}^{\mathcal{R}_{\partial}}\left(  \mathbb{R}^{d}\right)  $ and $\left(
u_{\varepsilon}\right)  _{\varepsilon}\in\mathcal{X}_{\mathcal{S}%
}^{\mathcal{R}_{\partial}}\left(  \mathbb{R}^{d}\right)  $ a representative of
$u$. There exists a sequence $N\in\mathcal{R}$ such that $\mu_{r,m}\left(
u_{\varepsilon}\right)  =\mathrm{O}\left(  \varepsilon^{-N(r)}\right)  $
$\mathrm{as}\;\varepsilon\rightarrow0$, for all $r\in\mathbb{N}$. According to
the stability of regular sets, there exists a sequence $N^{\prime}%
\in\mathcal{R}$ such that
\[
\forall l\in\mathbb{N}\ \ \ N(l+d+1)\leq N^{\prime}(l).
\]
Lemma \ref{GregLemFond} implies that $\mu_{q,l}\left(  \widehat{u}%
_{\varepsilon}\right)  =\mathrm{O}\left(  \varepsilon^{-N^{\prime}(q)}\right)
\;\mathrm{as}\;\varepsilon\rightarrow0$, for all $l\in\mathbb{N}$. Thus,
$\mathcal{F}\left(  u\right)  \in\mathcal{G}_{\mathcal{S}}^{\mathcal{R}_{u}%
}\left(  \mathbb{R}^{d}\right)  $.

So, we proved the inclusions of the sets in the left hand side of relations
(\ref{GregExcProp}), into the sets of the right hand side. The equalities
follow directly from a similar study with the inverse Fourier transform.
\end{proof}

\begin{example}
\label{Gregexchange}Take $\mathcal{R=}\mathbb{R}_{+}^{\mathbb{N}}$. We get
$\mathcal{F}\left(  \mathcal{G}_{\mathcal{S}}^{\mathcal{\partial}}\left(
\mathbb{R}^{d}\right)  \right)  =\mathcal{G}_{\mathcal{S}}^{u}\left(
\mathbb{R}^{d}\right)  $ and $\mathcal{F}\left(  \mathcal{G}_{\mathcal{S}}%
^{u}\left(  \mathbb{R}^{d}\right)  \right)  =\mathcal{G}_{\mathcal{S}%
}^{\mathcal{\partial}}\left(  \mathbb{R}^{d}\right)  $, result which is
closely related to the classical exchange theorem between $\mathcal{O}%
_{M}\left(  \mathbb{R}^{d}\right)  $ and $\mathcal{O}_{C}^{\prime}\left(
\mathbb{R}^{d}\right)  $.
\end{example}

Indeed, take $u\in\mathcal{O}_{C}^{\prime}\left(  \mathbb{R}^{d}\right)  $ and
consider $\left(  u_{\varepsilon}\right)  _{\varepsilon}=\left(  u\ast
\rho_{\varepsilon}\right)  _{\varepsilon}$ which is a representative of its
image by the embedding $\iota_{\mathcal{S}}$. Its Fourier image $\mathcal{F}%
\left(  \iota_{\mathcal{S}}\left(  u\right)  \right)  =\left[  \left(
\widehat{u}\widehat{\rho}_{\varepsilon}\right)  _{\varepsilon}\right]
_{\mathcal{S}}$ belongs to $\mathcal{G}_{\mathcal{S}}\left(  \mathbb{R}%
^{d}\right)  $, with $\widehat{u}\in\mathcal{O}_{M}\left(  \mathbb{R}%
^{d}\right)  $ and $\widehat{\rho}_{\varepsilon}\in\mathcal{S}\left(
\mathbb{R}^{d}\right)  $. As $\lim_{\varepsilon\rightarrow0}\widehat{\rho
}_{\varepsilon}=1$, we get $\lim_{\varepsilon\rightarrow0}\left(  \widehat
{u}\widehat{\rho}_{\varepsilon}\right)  _{\varepsilon}\in\mathcal{O}%
_{M}\left(  \mathbb{R}^{d}\right)  $. (For those limits, we consider
$\mathcal{O}_{M}\left(  \mathbb{R}^{d}\right)  $ equipped with its usual
topology: See \cite{Horv}, \cite{Schwartz1}.) This shows the consistency of
our result with the classical one. The generalized function $\mathcal{F}%
\left(  \iota_{\mathcal{S}}\left(  u\right)  \right)  $ belongs to a space of
rapidly decreasing generalized functions, but the limit of its representatives
when $\varepsilon\rightarrow0$ is in a space of functions of moderate growth.

\begin{corollary}
\label{GregFregCor}(\textbf{Regularity theorem}) We have: $\mathcal{F}\left(
\mathcal{G}_{\mathcal{S}}^{\mathcal{\infty}}\left(  \mathbb{R}^{d}\right)
\right)  =\mathcal{G}_{\mathcal{S}}^{\mathcal{\infty}}\left(  \mathbb{R}%
^{d}\right)  $.
\end{corollary}

\begin{proof}
Apply theorem \ref{GregExchangethm} with $\mathcal{R}=\mathcal{B}$, the set of
bounded sequences, for which $\mathcal{B}_{u}=\mathcal{B}_{\partial}$%
.\medskip\smallskip
\end{proof}

We can now complete diagram \ref{GregDiagramDU1} in the case of $\Omega
=\mathbb{R}^{d}$:%
\begin{equation}%
\begin{array}
[c]{ccccccc}
&  & \mathcal{G}_{\mathcal{S}}^{\mathcal{R}_{\partial}}\left(  \mathbb{R}%
^{d}\right)  & \longrightarrow & \mathcal{G}_{\mathcal{S}}^{\partial}\left(
\mathbb{R}^{d}\right)  &  & \\
\mathcal{G}_{\mathcal{S}}^{\infty}\left(  \mathbb{R}^{d}\right)  &
\begin{array}
[c]{c}%
\nearrow\\
\searrow
\end{array}
& \updownarrow\mathcal{F} &  & \updownarrow\mathcal{F} &
\begin{array}
[c]{c}%
\searrow\\
\nearrow
\end{array}
& \mathcal{G}_{\mathcal{S}}\left(  \mathbb{R}^{d}\right) \\
&  & \mathcal{G}_{\mathcal{S}}^{\mathcal{R}_{u}}\left(  \mathbb{R}^{d}\right)
& \longrightarrow & \mathcal{G}_{\mathcal{S}}^{u}\left(  \mathbb{R}^{d}\right)
&  &
\end{array}
. \label{GregDiagramDU2}%
\end{equation}

An interesting consequence of corollary \ref{GregFregCor} is the following
property, also proved in \cite{Garetto}, which is the equivalent for rapidly
decreasing generalized functions of the result mentioned in the introduction
for the $\mathcal{G}^{\infty}$ regularity ($\mathcal{D}^{\prime}\left(
\Omega\right)  \cap\mathcal{G}^{\infty}\left(  \Omega\right)  =\mathrm{C}%
^{\infty}\left(  \Omega\right)  $, \cite{Ober1}).

\begin{proposition}
\label{GSFLequalSta}We have $\mathcal{O}_{C}^{\prime}\left(  \mathbb{R}%
^{d}\right)  \cap\mathcal{G}_{\mathcal{S}}^{\infty}\left(  \mathbb{R}%
^{d}\right)  =\mathcal{S}\left(  \mathbb{R}^{d}\right)  $.
\end{proposition}

\begin{proof}
We follow here the ideas of \cite{NePiSc} for the proof of the above mentioned
result about $\mathcal{G}^{\infty}\left(  \mathbb{R}^{d}\right)  $. Let $u$ be
in $\mathcal{O}_{C}^{\prime}\left(  \mathbb{R}^{d}\right)  $ and set $\left(
u_{\varepsilon}\right)  _{\varepsilon}=\left(  u\ast\rho_{\varepsilon}\right)
_{\varepsilon}$. By assumption $\left[  \left(  u\ast\rho_{\varepsilon
}\right)  _{\varepsilon}\right]  _{\mathcal{S}}$ is in $\mathcal{G}%
_{\mathcal{S}}^{\infty}\left(  \mathbb{R}^{d}\right)  $. According to
corollary \ref{GregFregCor}, $\mathcal{F}_{S}\left(  \left[  \left(  u\ast
\rho_{\varepsilon}\right)  _{\varepsilon}\right]  _{\mathcal{S}}\right)  $ is
also in $\mathcal{G}_{\mathcal{S}}^{\infty}\left(  \mathbb{R}^{d}\right)  $.
It follows that there exists $N\in\mathbb{N}$ such that%
\[
\forall q\in\mathbb{N},\;\;\exists C_{q}>0\;\;\;\sup_{\xi\in\mathbb{R}^{d}%
}\left(  1+\left\vert \xi\right\vert \right)  ^{q}\left\vert \widehat
{u}\left(  \xi\right)  \widehat{\rho}_{\varepsilon}\left(  \xi\right)
\right\vert \leq C_{q}\varepsilon^{-N},\;\text{for }\varepsilon\text{ small
enough.}%
\]
By choice of $\rho$, $\widehat{\rho}_{\varepsilon}$ is an element of
$\mathcal{D}\left(  \mathbb{R}^{d}\right)  $.\ Moreover, a straightforward
calculation shows that $\widehat{\rho}_{\varepsilon}\left(  \xi\right)
=\widehat{\rho}\left(  \varepsilon\xi\right)  $, for all $\xi\in\mathbb{R}%
^{d}$, with $\widehat{\rho}$ equal to $1$ on a neighborhood of $0$. It follows
that, for all $q\in\mathbb{N}$, we have
\begin{align*}
\left.  \forall\xi\in\mathbb{R}^{d}\;\;\;\left(  1+\left\vert \xi\right\vert
\right)  ^{q}\left\vert \widehat{u}\left(  \xi\right)  \right\vert \right.
&  \leq\left(  1+\left\vert \xi\right\vert \right)  ^{q}\left\vert \widehat
{u}\left(  \xi\right)  \right\vert \left(  \left\vert 1-\widehat{\rho}\left(
\varepsilon\xi\right)  \right\vert +\left\vert \widehat{\rho}\left(
\varepsilon\xi\right)  \right\vert \right) \\
&  \leq\left(  1+\left\vert \xi\right\vert \right)  ^{q}\left\vert \widehat
{u}\left(  \xi\right)  \right\vert \left\vert 1-\widehat{\rho}\left(
\varepsilon\xi\right)  \right\vert +C_{q}\varepsilon^{-N}.
\end{align*}
Since $1-\widehat{\rho}\left(  \varepsilon\xi\right)  =\widehat{\rho}\left(
0\right)  -\widehat{\rho}\left(  \varepsilon\xi\right)  =-\varepsilon\xi
\int_{0}^{1}\widehat{\rho}^{\prime}\left(  \varepsilon\xi t\right)
\,\mathrm{d}t$, with $\widehat{\rho}^{\prime}$ bounded, there exists a
constant $C>0$ such that
\[
\left.  \forall\xi\in\mathbb{R}^{d}\;\;\;\left(  1+\left\vert \xi\right\vert
\right)  ^{q}\left\vert \widehat{u}\left(  \xi\right)  \right\vert \right.
\leq C\left(  1+\left\vert \xi\right\vert \right)  ^{q}\left\vert \widehat
{u}\left(  \xi\right)  \right\vert \varepsilon\left\vert \xi\right\vert
+C_{q}\varepsilon^{-N}.
\]
As $\widehat{u}$ is in $\mathcal{O}_{M}\left(  \mathbb{R}^{d}\right)  $, there
exist $m\in\mathbb{N}$ and a constant $C_{1}>0$ such that%
\[
\sup\nolimits_{\xi\in\mathbb{R}^{d}}\left(  1+\left\vert \xi\right\vert
\right)  ^{-m+1}\left\vert \widehat{u}\left(  \xi\right)  \right\vert \leq
C_{1}.
\]
Therefore, by setting $C_{2}=\max\left(  CC_{1},C_{q}\right)  $, we get%
\begin{align*}
\left.  \forall\xi\in\mathbb{R}^{d}\;\;\;\left(  1+\left\vert \xi\right\vert
\right)  ^{q}\left\vert \widehat{u}\left(  \xi\right)  \right\vert \right.
&  \leq C_{2}\left(  \left(  1+\left\vert \xi\right\vert \right)
^{q+m-1}\varepsilon\left\vert \xi\right\vert +\varepsilon^{-N}\right) \\
&  \leq C_{2}\left(  \left(  1+\left\vert \xi\right\vert \right)
^{q+m}\varepsilon+\varepsilon^{-N}\right)  .
\end{align*}
By minimizing the function $f_{\xi}:\varepsilon\mapsto\left(  1+\left\vert
\xi\right\vert \right)  ^{q+m}\varepsilon+\varepsilon^{-N}$, we get the
existence of a constant $C_{3}>0$ such that%
\[
\left.  \forall\xi\in\mathbb{R}^{d}\;\;\;\left(  1+\left\vert \xi\right\vert
\right)  ^{q}\left\vert \widehat{u}\left(  \xi\right)  \right\vert \right.
\leq C_{3}\left(  \left(  1+\left\vert \xi\right\vert \right)  ^{N\left(
q+m\right)  /(N+1)}\right)  ,
\]
and
\[
\left.  \forall\xi\in\mathbb{R}^{d}\;\;\;\left\vert \widehat{u}\left(
\xi\right)  \right\vert \right.  \leq C_{3}\left(  \left(  1+\left\vert
\xi\right\vert \right)  ^{-q/(N+1)+mN/(N+1)}\right)  ,
\]
for all $q\in\mathbb{N}$. ($m$ only depends on $u$.) Treating the derivatives
in the same way, we obtain the same type of estimates. Therefore $\widehat{u}$
and its derivatives are rapidly decreasing.\ This shows that $\mathcal{O}%
_{C}^{\prime}\left(  \Omega\right)  \cap\mathcal{G}_{\mathcal{S}}^{\infty
}\left(  \Omega\right)  \subset\mathcal{S}\left(  \Omega\right)  $.\ As the
other inclusion is obvious, our claim is proved.
\end{proof}

\section{Global regularity of compactly supported generalized
functions\label{Glimsec5}}

\subsection{$\mathrm{C}^{\infty}$-regularity for compactly supported
distributions}

In order to render easier the comparison between the distributional case and
the generalized case, we are going to recall the classical theorem and
complete it by some equivalent statements.

\begin{theorem}
\label{RDGFCinfreg}For $u$ in $\mathcal{E}^{\prime}\left(  \mathbb{R}%
^{d}\right)  $, the following equivalences hold:%
\[%
\begin{tabular}
[c]{lll}%
$\left(  i\right)  ~\,u\in\mathrm{C}^{\infty}\left(  \mathbb{R}^{d}\right)  $
& $\Leftrightarrow\left(  ii\right)  $ & $\!\!\mathcal{F}\left(  u\right)
\in\mathcal{S}\left(  \mathbb{R}^{d}\right)  $\\
& $\Leftrightarrow\left(  iii\right)  $ & $\!\!\mathcal{F}\left(  u\right)
\in\mathcal{O}_{C}^{\prime}\left(  \mathbb{R}^{d}\right)  $\\
& $\Leftrightarrow\left(  iv\right)  $ & $\!\!\mathcal{F}\left(  u\right)
\in\mathcal{O}_{M}^{\prime}\left(  \mathbb{R}^{d}\right)  $\\
& $\Leftrightarrow\left(  v\right)  $ & $\!\!\mathcal{F}\left(  u\right)
\in\mathcal{O}_{C}^{\prime}\left(  \mathbb{R}^{d}\right)  .$%
\end{tabular}
\]

\end{theorem}

\begin{proof}
Equivalence $\left(  i\right)  \Leftrightarrow\left(  ii\right)  $ is the
classical result.\ The trivial inclusion $\mathcal{S}\left(  \mathbb{R}%
^{d}\right)  \subset\mathcal{S}_{\ast}\left(  \mathbb{R}^{d}\right)  $ shows
$\left(  ii\right)  \Rightarrow\left(  iii\right)  $. Then, the structure of
elements of $\mathcal{O}_{M}^{\prime}\left(  \mathbb{R}^{d}\right)  $
\cite{Ortner} shows that $\mathcal{S}_{\ast}\left(  \Omega\right)  $ is
canonically embedded in $\mathcal{O}_{M}^{\prime}\left(  \mathbb{R}%
^{d}\right)  $: This shows $\left(  iii\right)  \Rightarrow\left(  iv\right)
$.\ As $\mathcal{O}_{M}^{\prime}\left(  \mathbb{R}^{d}\right)  \subset
\mathcal{O}_{C}^{\prime}\left(  \mathbb{R}^{d}\right)  $, $\left(  iv\right)
\Rightarrow\left(  v\right)  $ is obvious. For $\left(  v\right)
\Rightarrow\left(  i\right)  $, note that $\mathcal{F}\left(  u\right)  $
belongs to $\mathcal{O}_{M}\left(  \mathbb{R}^{d}\right)  $ and better to
$\mathcal{O}_{C}\left(  \mathbb{R}^{d}\right)  $ since $u$ is in
$\mathcal{E}^{\prime}\left(  \mathbb{R}^{d}\right)  $. (This last assertion is
a refinement of the classical previous one.) Then, if $\left(  v\right)  $
holds, $\mathcal{F}\left(  u\right)  $ is in $\mathcal{O}_{C}\left(
\mathbb{R}^{d}\right)  \cap\mathcal{O}_{C}^{\prime}\left(  \mathbb{R}%
^{d}\right)  $ which is equal to $\mathcal{S}\left(  \mathbb{R}^{d}\right)  $
\cite{Ortner}. Then $\left(  ii\right)  $ holds.\medskip
\end{proof}

Theorem \ref{RDGFCinfreg} shows, at least, that there is no need to consider
spaces of functions with all the derivatives rapidly decreasing to
characterize elements of $\mathcal{E}^{\prime}\left(  \mathbb{R}^{d}\right)  $
which are $\mathrm{C}^{\infty}$. In fact, we can only consider functions
rapidly decreasing, with no other hypothesis on the derivatives. A similar
situation holds for generalized functions, justifying the introduction of
rough generalized functions in the following subsection.

\subsection{Rough rapidly decreasing generalized functions}

\subsubsection{Definitions}

Let $\mathcal{R}$ be a regular subset of $\mathbb{R}_{+}^{\mathbb{N}}$ and
$\Omega$ an open subset or $\mathbb{R}^{d}$.\ Set%
\begin{align}
\mathcal{S}_{\ast}\left(  \Omega\right)   &  =\left\{  f\in\mathrm{C}^{\infty
}\left(  \Omega\right)  \,\left\vert \,\forall q\in\mathbb{N}\;\;\mu
_{q,0}\left(  f\right)  <+\infty\right.  \right\}  ,\nonumber\\
\mathcal{X}_{\mathcal{S}_{\ast}}^{\mathcal{R}}\left(  \Omega\right)   &
=\left\{  \left(  f_{\varepsilon}\right)  _{\varepsilon}\in\mathcal{S}_{\ast
}\left(  \Omega\right)  ^{\left(  0,1\right]  }\,\left\vert \,\exists
N\in\mathcal{R},\ \forall q\in\mathbb{N}\;\;\mu_{q,0}\left(  f_{\varepsilon
}\right)  =\mathrm{O}\left(  \varepsilon^{-N(q)}\right)  \;\mathrm{as}%
\;\varepsilon\rightarrow0\right.  \right\}  ,\label{GlimXiSDef2}\\
\mathcal{N}_{\mathcal{S}_{\ast}}\left(  \Omega\right)   &  =\left\{  \left(
f_{\varepsilon}\right)  _{\varepsilon}\in\mathcal{S}_{\ast}\left(
\Omega\right)  ^{\left(  0,1\right]  }\,\left\vert \,\forall N\in
\mathbb{R}_{+}^{\mathbb{N}},\;\forall l\in\mathbb{N}\;\;\mu_{q,0}\left(
f_{\varepsilon}\right)  =\mathrm{O}\left(  \varepsilon^{N(l)}\right)
\;\mathrm{as}\;\varepsilon\rightarrow0\right.  \right\}  .\nonumber
\end{align}
One can show that $\mathcal{X}_{\mathcal{S}_{\ast}}^{\mathcal{R}}\left(
\Omega\right)  $ is a subalgebra of $\mathcal{S}_{\ast}\left(  \Omega\right)
^{\left(  0,1\right]  }$ and $\mathcal{N}_{\mathcal{S}_{\ast}}\left(
\Omega\right)  $ an ideal of $\mathcal{X}_{\mathcal{S}_{\ast}}^{\mathcal{R}%
}\left(  \Omega\right)  $. (In fact, these spaces fit in the general scheme of
construction of Colombeau type algebra \cite{AntRad}, \cite{DHPV1},
\cite{JAM2}.)

\begin{definition}
\label{GlimdefRg1}The space $\mathcal{G}_{\mathcal{S}_{\ast}}^{\mathcal{R}%
}\left(  \Omega\right)  =\mathcal{X}_{\mathcal{S}_{\ast}}^{\mathcal{R}}\left(
\Omega\right)  /\mathcal{N}_{\mathcal{S}_{\ast}}\left(  \Omega\right)  $ is
called the algebra of $\mathcal{R}$\emph{-regular\ rough rapidly decreasing
generalized functions.}
\end{definition}

\begin{example}
\label{GregRoughExamp1}Taking $\mathcal{R}=\mathbb{R}_{+}^{\mathbb{N}}$, we
obtain the space $\mathcal{G}_{\mathcal{S}_{\ast}}\left(  \Omega\right)  $ of
\emph{rough rapidly decreasing generalized functions}.
\end{example}

\begin{example}
\label{GregRoughExamp2}Taking $\mathcal{R}=\mathcal{B}$, the set of bounded
sequences, we obtain the space $\mathcal{G}_{\mathcal{S}_{\ast}}^{\infty
}\left(  \Omega\right)  $, of \emph{regular rough rapidly decreasing
generalized functions}.
\end{example}

Lemma \ref{GregLemmaBox} implies immediately the following proposition:

\begin{proposition}
\label{GregGSRGSRstar}If the open set $\Omega$ is a box and $\mathcal{R}%
^{\prime}$ a regular subset of $\mathbb{R}_{+}^{\mathbb{N}^{2}}$, then
$\mathcal{G}_{\mathcal{S}}^{\mathcal{R}^{\prime}}\left(  \Omega\right)  $ is
included in $\mathcal{G}_{\mathcal{S}_{\ast}}^{\mathcal{R}_{0}^{\prime}%
}\left(  \Omega\right)  $, where $\mathcal{R}_{0}^{\prime}$ is the (regular)
subset of $\mathbb{R}_{+}^{\mathbb{N}}$ equal to $\left\{  N\left(
\cdot,0\right)  \text{, }N\in\mathcal{R}^{\prime}\right\}  $.
\end{proposition}

\begin{example}
\label{GregGSDGSinf}If $\Omega$ is a box, for all $\mathcal{R}\subset
\mathbb{R}_{+}^{\mathbb{N}}$, $\mathcal{G}_{\mathcal{S}}^{\mathcal{R}%
_{\partial}}\left(  \Omega\right)  $ is included in $\mathcal{G}%
_{\mathcal{S}_{\ast}}^{\mathcal{R}}\left(  \Omega\right)  .$
\end{example}

Indeed, $\mathcal{R}_{\partial}=\mathcal{R}\otimes\left\{  1\right\}  $, which
implies that $\left(  \mathcal{R}_{\partial}\right)  _{0}=\mathcal{R}$. Let us
quote two other examples of application of proposition \ref{GregGSRGSRstar}.

\begin{corollary}
\label{GregRoughExamp1bis}If the open set $\Omega$ is a box, then
\newline$\left(  i\right)  $~$\mathcal{G}_{\mathcal{S}}\left(  \Omega\right)
$, obtained for $\mathcal{R}^{\prime}=\mathbb{R}_{+}^{\mathbb{N}^{2}}$, is
included in $\mathcal{G}_{\mathcal{S}_{\ast}}\left(  \Omega\right)  $%
.\newline$\left(  ii\right)  $~$\mathcal{G}_{\mathcal{S}}^{\infty}\left(
\Omega\right)  $, obtained for $\mathcal{R}^{\prime}=\mathcal{B}^{\prime}$, is
included in $\mathcal{G}_{\mathcal{S}_{\ast}}^{\infty}\left(  \Omega\right)  $.
\end{corollary}

Indeed, $\left(  i\right)  $ (\textit{resp}.\ $\left(  ii\right)  $) holds,
since $\left(  \mathbb{R}_{+}^{\mathbb{N}^{2}}\right)  _{0}=\mathbb{R}%
_{+}^{\mathbb{N}}$ (\textit{resp}. $\left(  \mathcal{B}^{\prime}\right)
_{0}=\mathcal{B}^{\prime}$). Note that the proof of proposition
\ref{GSFLequalSta} shows that%
\[
\mathcal{G}_{\mathcal{S}_{\ast}}^{\infty}\left(  \mathbb{R}^{d}\right)
\cap\mathcal{O}_{C}^{\prime}\left(  \mathbb{R}^{d}\right)  =\mathcal{S}_{\ast
}\left(  \mathbb{R}^{d}\right)  .
\]

We turn to the question of embeddings. First, the structure of elements of
$\mathcal{O}_{C}^{\prime}\left(  \mathbb{R}^{d}\right)  $ (\cite{Horv},
\cite{Ortner}, \cite{Schwartz1}) shows that $\mathcal{S}_{\ast}\left(
\mathbb{R}^{d}\right)  $ is canonically embedded in $\mathcal{O}_{C}^{\prime
}\left(  \mathbb{R}^{d}\right)  $.\ The embedding of $\mathcal{S}\left(
\mathbb{R}^{d}\right)  $ into $\mathcal{G}_{\mathcal{S}_{\ast}}\left(
\mathbb{R}^{d}\right)  $ is done by the canonical injective map
\[
\sigma_{\mathcal{S}_{\ast}}:\mathcal{S}_{\ast}\left(  \mathbb{R}^{d}\right)
\rightarrow\mathcal{G}_{\mathcal{S}_{\ast}}\left(  \mathbb{R}^{d}\right)
\;\;\ f\mapsto\left(  f_{\varepsilon}\right)  _{\varepsilon}+\mathcal{N}%
_{\mathcal{S}_{\ast}}\left(  \mathbb{R}^{d}\right)  \ \ \mathrm{with\;}%
f_{\varepsilon}=f\text{ for }\varepsilon\in\left(  0,1\right]  .
\]
Finally, a simplification of the proofs of theorems \ref{CSTGThemb},
\ref{CSTGThembcom} and proposition \ref{GregImbedOPM} leads to the following
theorem, where $\left(  \rho_{\varepsilon}\right)  _{\varepsilon}$ is defined
by (\ref{CSTGSBinfini}) and (\ref{CSTGSBinfini2}).

\begin{theorem}
\label{RGDFembrough}~\newline$\left(  i\right)  $~The map
\[
\iota_{\mathcal{S}_{\ast}}:\mathcal{O}_{C}^{\prime}\left(  \mathbb{R}%
^{d}\right)  \rightarrow\mathcal{G}_{\mathcal{S}_{\ast}}\left(  \mathbb{R}%
^{d}\right)  \;\;\;\;\;\;u\mapsto\left(  u\ast\rho_{\varepsilon}\right)
_{\varepsilon}+\mathcal{N}_{\mathcal{S}_{\ast}}\left(  \mathbb{R}^{d}\right)
\]
is a linear embedding which commutes with partial derivatives.\newline$\left(
ii\right)  $~We have: $\iota_{\mathcal{S}_{\ast}\left\vert \mathcal{S}_{\ast
}\left(  \mathbb{R}^{d}\right)  \right.  }=\sigma_{\mathcal{S}_{\ast}}%
.$\newline$\left(  iii\right)  $~We have: $\iota_{\mathcal{S}_{\ast}}\left(
\mathcal{O}_{M}^{\prime}\left(  \mathbb{R}^{d}\right)  \right)  \subset
\mathcal{G}_{\mathcal{S}_{\ast}}^{\infty}\left(  \mathbb{R}^{d}\right)  $.
\end{theorem}

\begin{remark}
\label{RGDFRemarkcom}Theorems \ref{CSTGThemb}, \ref{CSTGThembcom} and
\ref{RGDFembrough} combined together show that all the arrows are injective
and all diagrams commutative in the following schemes:
\[%
\begin{tabular}
[c]{ccccc}%
$\mathcal{S}\left(  \mathbb{R}^{d}\right)  $ &  & $\longrightarrow$ &  &
$\mathcal{S}_{\ast}\left(  \mathbb{R}^{d}\right)  $\\
& $\searrow$ &  & $\swarrow$ & \\
$\downarrow$ &  & $\mathcal{O}_{C}^{\prime}\left(  \mathbb{R}^{d}\right)  $ &
& $\downarrow$\\
& $\swarrow$ &  & $\searrow$ & \\
$\mathcal{G}_{\mathcal{S}}\left(  \mathbb{R}^{d}\right)  $ &  &
$\longrightarrow$ &  & $\mathcal{G}_{\mathcal{S}_{\ast}}\left(  \mathbb{R}%
^{d}\right)  $%
\end{tabular}
\ \ \ \ \ \ \ \ \ \ \ \ \ \ \
\begin{tabular}
[c]{ccccc}%
$\mathcal{S}\left(  \mathbb{R}^{d}\right)  $ &  & $\longrightarrow$ &  &
$\mathcal{S}_{\ast}\left(  \mathbb{R}^{d}\right)  $\\
& $\searrow$ &  & $\swarrow$ & \\
$\downarrow$ &  & $\mathcal{O}_{M}^{\prime}\left(  \mathbb{R}^{d}\right)  $ &
& $\downarrow$\\
& $\swarrow$ &  & $\searrow$ & \\
$\mathcal{G}_{\mathcal{S}}^{u}\left(  \mathbb{R}^{d}\right)  $ &  &
$\longrightarrow$ &  & $\mathcal{G}_{\mathcal{S}_{\ast}}^{\infty}\left(
\mathbb{R}^{d}\right)  $%
\end{tabular}
\ .
\]

\end{remark}

\subsubsection{Fourier transform in $\mathcal{G}_{\mathcal{S}_{\ast}}\left(
\mathbb{R}^{d}\right)  $}

We need in the sequel to define a Fourier transform (or an inverse Fourier
transform) in $\mathcal{G}_{\mathcal{S}_{\ast}}^{\mathcal{R}}\left(
\mathbb{R}^{d}\right)  $. This is done in the following way.\ Set, for any
regular subspace $\mathcal{R}$ of $\mathbb{R}_{+}^{\mathbb{N}}$,
\begin{align*}
\mathcal{X}_{\mathcal{B}}\left(  \Omega\right)   &  =\left\{  \left(
f_{\varepsilon}\right)  _{\varepsilon}\in\mathrm{C}^{\infty}\left(
\Omega\right)  ^{\left(  0,1\right]  }\,\left\vert \,\exists N\in
\mathbb{R}_{+}^{\mathbb{N}},\;\forall l\in\mathbb{N}\;\;\mu_{0,l}\left(
f_{\varepsilon}\right)  =\mathrm{O}\left(  \varepsilon^{-N(l)}\right)
\;\mathrm{as}\;\varepsilon\rightarrow0\right.  \right\}  ,\\
\mathcal{X}_{\mathcal{B}}^{\mathcal{R}}\left(  \Omega\right)   &  =\left\{
\left(  f_{\varepsilon}\right)  _{\varepsilon}\in\mathrm{C}^{\infty}\left(
\Omega\right)  ^{\left(  0,1\right]  }\,\left\vert \,\exists N\in
\mathcal{R},\;\forall l\in\mathbb{N}\;\;\mu_{0,l}\left(  f_{\varepsilon
}\right)  =\mathrm{O}\left(  \varepsilon^{-N(l)}\right)  \;\mathrm{as}%
\;\varepsilon\rightarrow0\right.  \right\}  ,\\
\mathcal{N}_{\mathcal{B}}\left(  \Omega\right)   &  =\left\{  \left(
f_{\varepsilon}\right)  _{\varepsilon}\in\mathrm{C}^{\infty}\left(
\Omega\right)  ^{\left(  0,1\right]  }\,\left\vert \,\forall N\in
\mathbb{R}_{+}^{\mathbb{N}},\;\forall l\in\mathbb{N}\;\;\mu_{0,l}\left(
f_{\varepsilon}\right)  =\mathrm{O}\left(  \varepsilon^{N(l)}\right)
\;\mathrm{as}\;\varepsilon\rightarrow0\right.  \right\}  .
\end{align*}

According to the general scheme of construction of Colombeau type algebras,
$\mathcal{G}_{\mathcal{B}}\left(  \Omega\right)  =\mathcal{X}_{\mathcal{B}%
}\left(  \Omega\right)  /\mathcal{N}_{\mathcal{B}}\left(  \Omega\right)  $ is
an algebra, named the \emph{algebra of bounded generalized functions}.
Moreover, $\mathcal{X}_{\mathcal{B}}^{\mathcal{R}}\left(  \Omega\right)  $ is
a subalgebra of $\mathcal{X}_{\mathcal{B}}\left(  \Omega\right)  $. (The proof
is similar to the one of proposition \ref{GlimPropFReg}.) The space
$\mathcal{G}_{\mathcal{B}}^{\mathcal{R}}\left(  \Omega\right)  =\mathcal{X}%
_{\mathcal{B}}^{\mathcal{R}}\left(  \Omega\right)  /\mathcal{N}_{\mathcal{B}%
}\left(  \Omega\right)  $ is called the space of $\mathcal{R}$-\emph{regular
bounded generalized functions}.

\begin{notation}
We shall note $\left[  \left(  f_{\varepsilon}\right)  _{\varepsilon}\right]
_{\mathcal{B}}$ the class of $\left(  f_{\varepsilon}\right)  _{\varepsilon}$
in $\mathcal{G}_{\mathcal{B}}^{\mathcal{R}}\left(  \Omega\right)  $.
\end{notation}

\begin{remark}
\label{GregGCembGBR}One can verify that $\mathcal{G}_{C}\left(  \Omega\right)
$ (resp. $\mathcal{G}_{C}^{\mathcal{R}}\left(  \Omega\right)  $) is embedded
in $\mathcal{G}_{\mathcal{B}}\left(  \Omega\right)  $ (resp. $\mathcal{G}%
_{\mathcal{B}}^{\mathcal{R}}\left(  \Omega\right)  $).
\end{remark}

\begin{proposition}
\label{GlimLmTF2}$~$\newline$\left(  i\right)  ~$For all $u\in\mathcal{G}%
_{\mathcal{S}_{\ast}}\left(  \mathbb{R}^{d}\right)  $ and $\left(
u_{\varepsilon}\right)  _{\varepsilon}\in\mathcal{X}_{\mathcal{S}_{\ast}%
}\left(  \mathbb{R}^{d}\right)  $ a representative of $u$, the expression
\begin{equation}
\widehat{u}:\left[  \widehat{u}_{\varepsilon}=\left(  \xi\mapsto\int
\mathrm{e}^{-ix\xi}u_{\varepsilon}\left(  x\right)  \,\mathrm{d}x\right)
_{\varepsilon}\right]  _{\mathcal{B}} \label{GlimEq8}%
\end{equation}
defines an element of $\mathcal{G}_{\mathcal{B}}\left(  \Omega\right)  $
depending only on $u$.\newline$\left(  ii\right)  ~$For any regular subspace
$\mathcal{R}$ of $\mathbb{R}_{+}^{\mathbb{N}}$ and $\left(  u_{\varepsilon
}\right)  _{\varepsilon}\in\mathcal{X}_{\mathcal{S}}^{\mathcal{R}}\left(
\Omega\right)  $, we have $\left(  \widehat{u}_{\varepsilon}\right)
_{\varepsilon}\in\mathcal{X}_{\mathcal{B}}^{\mathcal{R}}\left(  \Omega\right)
$.
\end{proposition}

The \textbf{proof} of proposition \ref{GlimLmTF2} is mainly a consequence of
lemma \ref{GregLemFond}.

\noindent\textit{Assertion} $\left(  i\right)  .~$Take $u\in\mathcal{G}%
_{\mathcal{S}_{\ast}}\left(  \mathbb{R}^{d}\right)  $ and $\left(
u_{\varepsilon}\right)  _{\varepsilon}\in\mathcal{X}_{\mathcal{S}_{\ast}%
}\left(  \mathbb{R}^{d}\right)  $ a representative of $u$. Then lemma
\ref{GregLemFond} (applied with $q=0$) implies that
\begin{equation}
\forall l\in\mathbb{N},\ \exists C_{l}>0,\ \forall\varepsilon\in\left(
0,1\right]  \ \ \ \ \mu_{0,l}\left(  \widehat{u}_{\varepsilon}\right)  \leq
C_{\alpha}\,\mu_{l+d+1,0}\left(  u_{\varepsilon}\right)  . \label{GlimEq9b}%
\end{equation}
This estimate shows that $\left(  \widehat{u}_{\varepsilon}\right)
_{\varepsilon}\in\mathcal{X}_{\mathcal{B}}\left(  \mathbb{R}^{d}\right)  $.
Indeed, if $\left(  u_{\varepsilon}\right)  _{\varepsilon}$ is in
$\mathcal{X}_{\mathcal{S}_{\ast}}\left(  \mathbb{R}^{d}\right)  $, there
exists a sequence $N\in\mathcal{R}$ such that $\mu_{l,0}\left(  u_{\varepsilon
}\right)  =\mathrm{O}\left(  \varepsilon^{-N(l)}\right)  \;\mathrm{as}%
\;\varepsilon\rightarrow0$ and setting $N^{\prime}:l\mapsto N\left(
l+d+1\right)  $, we get that $\mu_{0,l}\left(  \widehat{u}_{\varepsilon
}\right)  =\left(  \varepsilon^{-N^{\prime}(l)}\right)  \;\mathrm{as}%
\;\varepsilon\rightarrow0$. According to the overstability by translation of
the subset $\mathcal{R}$, $\left(  \widehat{u}_{\varepsilon}\right)
_{\varepsilon}$ belongs to $\mathcal{X}_{\mathcal{B}}\left(  \mathbb{R}%
^{d}\right)  $. Similar arguments show that, if $\left(  \eta_{\varepsilon
}\right)  _{\varepsilon}\in\mathcal{N}_{\mathcal{S}_{\ast}}\left(
\mathbb{R}^{d}\right)  $, then $\left(  \widehat{\eta}_{\varepsilon}\right)
_{\varepsilon}\in\mathcal{N}_{\mathcal{B}}\left(  \Omega\right)  $. Therefore,
relation (\ref{GlimEq8}) defines an element of $\mathcal{G}_{\mathcal{B}%
}\left(  \mathbb{R}^{d}\right)  $, depending only on $u$.\smallskip

\noindent\textit{Assertion} $\left(  ii\right)  $. The estimate
(\ref{GlimEq9b}) implies that the regularity of the sequences in the
definition of moderate elements transfers by Fourier transform from the space
index $q$ in the $\mathcal{S}_{\ast}$-type spaces to the derivative index $l$
in the Colombeau type space (here of bounded functions), showing our
claim.\smallskip\medskip

We define the\textbf{ Fourier transform} $\mathcal{F}_{\ast}$ on
$\mathcal{G}_{\mathcal{S}_{\ast}}\left(  \mathbb{R}^{d}\right)  \mathcal{\ }%
$by the formula%
\[
\mathcal{F}_{\ast}:\mathcal{G}_{\mathcal{S}_{\ast}}\left(  \mathbb{R}%
^{d}\right)  \rightarrow\mathcal{G}_{B}\left(  \mathbb{R}^{d}\right)
\;\;\;u\mapsto\left[  \left(  x\mapsto\int\mathrm{e}^{-ix\xi}u_{\varepsilon
}\left(  \xi\right)  \,\mathrm{d}\xi\right)  _{\varepsilon}\right]
_{\mathcal{B}}%
\]
where $\left(  u_{\varepsilon}\right)  _{\varepsilon}\in\mathcal{X}%
_{\mathcal{S}_{\ast}}\left(  \mathbb{R}^{d}\right)  $ is any representative of
$u$. (The inverse Fourier on $\mathcal{G}_{\mathcal{S}_{\ast}}\left(
\mathbb{R}^{d}\right)  $ is defined analogously.)\smallskip

The assertion $\left(  ii\right)  $ of proposition \ref{GlimLmTF2} implies:

\begin{proposition}
\label{GlimExthm1}(Small exchange theorem) We have: $\mathcal{F}\left(
\mathcal{G}_{\mathcal{S}_{\ast}}^{\mathcal{R}}\right)  \subset\mathcal{G}%
_{B}^{\mathcal{R}}\left(  \mathbb{R}^{d}\right)  $.
\end{proposition}

\subsection{$\mathcal{G}^{\mathcal{R}}$-regularity for compactly supported
generalized functions}

We have now all the elements to formulate and prove the following fundamental theorem:

\begin{theorem}
\label{GregGlobalTh1}Let $\mathcal{R}$ be regular subspace of $\mathbb{R}%
_{+}^{\mathbb{N}}$. For $u$ in $\mathcal{G}_{C}\left(  \mathbb{R}^{d}\right)
$, the following equivalences hold:%
\[%
\begin{tabular}
[c]{lll}%
$\left(  i\right)  ~u\in\mathcal{G}^{\mathcal{R}}\left(  \mathbb{R}%
^{d}\right)  $ & $\Leftrightarrow\left(  ii\right)  $ & $\!\!\mathcal{F}%
\left(  u\right)  \in\mathcal{G}_{\mathcal{S}}^{\mathcal{R}_{\partial}}\left(
\mathbb{R}^{d}\right)  \smallskip$\\
& $\Leftrightarrow\left(  iii\right)  $ & $\!\!\mathcal{F}\left(  u\right)
\in\mathcal{G}_{\mathcal{S}_{\ast}}^{\mathcal{R}}\left(  \mathbb{R}%
^{d}\right)  .$%
\end{tabular}
\]

\end{theorem}

\begin{proof}
~\newline$\left(  i\right)  \Rightarrow\left(  ii\right)  $~As $u$ is in
$\mathcal{G}_{C}\left(  \mathbb{R}^{d}\right)  \cap\mathcal{G}^{\mathcal{R}%
}\left(  \mathbb{R}^{d}\right)  =\mathcal{G}_{C}^{\mathcal{R}}\left(
\mathbb{R}^{d}\right)  $, $u$ is in $\mathcal{G}_{\mathcal{S}}^{\mathcal{R}%
_{u}}\left(  \mathbb{R}^{d}\right)  $ according to proposition
\ref{GregEmbGCRGSRU}. Then, applying theorem \ref{GregExchangethm},
$\mathcal{F}\left(  u\right)  $ is in $\mathcal{G}_{\mathcal{S}}%
^{\mathcal{R}_{\partial}}\left(  \mathbb{R}^{d}\right)  $.\smallskip

\noindent$\left(  ii\right)  \Rightarrow\left(  iii\right)  $~We have
$\mathcal{G}_{\mathcal{S}}^{\mathcal{R}_{\partial}}\left(  \mathbb{R}%
^{d}\right)  \subset\mathcal{G}_{\mathcal{S}_{\ast}}^{\mathcal{R}}\left(
\mathbb{R}^{d}\right)  $, according to example \ref{GregGSDGSinf}.\smallskip

\noindent$\left(  iii\right)  \Rightarrow\left(  i\right)  $~Let $u$ be in
$\mathcal{G}_{C}\left(  \mathbb{R}^{d}\right)  ,$ $\left(  u_{\varepsilon
}\right)  _{\varepsilon}$ be a representative of $u$ and $K$ a compact set
such that $\operatorname*{supp}u_{\varepsilon}\subset K$, for all
$\varepsilon$ in $\left(  0,1\right]  $. We have $\mathcal{F}_{\mathcal{S}%
}\left(  u\right)  =\left[  \left(  \widehat{u}_{\varepsilon}\right)
_{\varepsilon}\right]  _{\mathcal{G}_{\mathcal{S}}}$ where $\widehat{~}$
denotes the classical Fourier transform in $\mathcal{S}$. By assumption
$\mathcal{F}_{\mathcal{S}}\left(  u\right)  $ is in $\mathcal{G}%
_{\mathcal{S}_{\ast}}^{\mathcal{R}}\left(  \mathbb{R}^{d}\right)  $ and we can
consider its inverse Fourier transform $\mathcal{F}_{\ast}^{-1}$, with
$\mathcal{F}_{\ast}^{-1}\left(  \mathcal{F}_{\mathcal{S}}\left(  u\right)
\right)  $ in $\mathcal{G}_{B}^{\mathcal{R}}\left(  \mathbb{R}^{d}\right)  $
and
\[
\mathcal{F}_{\ast}^{-1}\left(  \mathcal{F}_{\mathcal{S}}\left(  u\right)
\right)  =\left[  \left(  \mathcal{F}^{-1}\left(  \widehat{u}_{\varepsilon
}\right)  \right)  _{\varepsilon}\right]  _{\mathcal{B}}.
\]
Using the classical isomorphism theorem in $\mathcal{S}$, we have
$\mathcal{F}^{-1}\left(  \widehat{u}_{\varepsilon}\right)  =u_{\varepsilon}$
for all $\varepsilon$ in $\left(  0,1\right]  $. Then
\[
\mathcal{F}_{\ast}^{-1}\left(  \mathcal{F}_{\mathcal{S}}\left(  u\right)
\right)  =\left[  \left(  u_{\varepsilon}\right)  _{\varepsilon}\right]
_{\mathcal{B}}.
\]
Since all the $u_{\varepsilon}$ have their support included in the same
compact set, we obviously have $\left[  \left(  u_{\varepsilon}\right)
_{\varepsilon}\right]  _{\mathcal{B}}=\iota_{C,B}\left(  u\right)  $ where
$\iota_{C,B}$ is the canonical embedding of $\mathcal{G}_{C}\left(
\mathbb{R}^{d}\right)  $ in $\mathcal{G}_{B}\left(  \mathbb{R}^{d}\right)  $.
Therefore, $u\in\mathcal{G}_{B}^{\mathcal{R}}\left(  \mathbb{R}^{d}\right)
\cap\mathcal{G}_{C}\left(  \mathbb{R}^{d}\right)  =\mathcal{G}^{\mathcal{R}%
}\left(  \mathbb{R}^{d}\right)  \cap\mathcal{G}_{C}\left(  \mathbb{R}%
^{d}\right)  $.
\end{proof}

\begin{example}
\label{GregGlobGinfR}The case $\mathcal{R}=\mathcal{B}$ in theorem
\ref{GregGlobalTh1} gives a characterization of the global $\mathcal{G}%
^{\infty}$-regularity of compactly supported generalized functions.
\end{example}

Moreover, we can refine theorem \ref{GregGlobalTh1} in this particular case
and prove:

\begin{theorem}
\label{RG-GlobRegGinf}For $u$ in $\mathcal{G}_{C}\left(  \mathbb{R}%
^{d}\right)  $, the following statements are equivalent:%
\[%
\begin{tabular}
[c]{lll}%
$\left(  i\right)  ~\,u\in\mathcal{G}^{\infty}\left(  \mathbb{R}^{d}\right)  $
& $\Leftrightarrow\left(  ii\right)  $ & $\!\!\mathcal{F}\left(  u\right)
\in\mathcal{G}_{\mathcal{S}}^{\infty}\left(  \mathbb{R}^{d}\right)
\smallskip$\\
& $\Leftrightarrow\left(  iii\right)  $ & $\!\!\mathcal{F}\left(  u\right)
\in\mathcal{G}_{\mathcal{S}}^{u}\left(  \mathbb{R}^{d}\right)  \smallskip$\\
& $\Leftrightarrow\left(  iv\right)  $ & $\!\!\mathcal{F}\left(  u\right)
\in\mathcal{G}_{\mathcal{S}_{\ast}}^{\infty}\left(  \mathbb{R}^{d}\right)  .$%
\end{tabular}
\]

\end{theorem}

Indeed, $\left(  i\right)  \Rightarrow\left(  ii\right)  $ and $\left(
iv\right)  \Rightarrow\left(  i\right)  $ follow directly from theorem
\ref{GregGlobalTh1} applied with $\mathcal{R}=\mathcal{B}$, since
$\mathcal{B}_{\partial}=\mathcal{B}^{\prime}$, the set of bounded elements of
$\mathbb{R}_{+}^{\mathbb{N}^{2}}$. For $\left(  ii\right)  \Rightarrow\left(
iii\right)  $, we have $\mathcal{G}_{\mathcal{S}}^{\infty}\left(
\mathbb{R}^{d}\right)  \subset\mathcal{G}_{\mathcal{S}}^{u}\left(
\mathbb{R}^{d}\right)  $.\ For $\left(  iii\right)  \Rightarrow\left(
iv\right)  $, we remark that $\mathcal{G}_{\mathcal{S}}^{u}\left(
\mathbb{R}^{d}\right)  $ is obtained with $\mathcal{R}^{\prime}=\left\{
1\right\}  \otimes\mathbb{R}_{+}^{\mathbb{N}}$ as regular subset of
$\mathbb{R}_{+}^{\mathbb{N}^{2}}$. This implies that $\left(  \mathcal{R}%
^{\prime}\right)  _{0}=\mathcal{B}^{\prime}$, with the notations of
proposition \ref{GregGSRGSRstar}.

\section{Local and microlocal $\mathcal{R}$-regularity\label{Glimsec6}}

We follow here the presentation of \cite{HorPDOT1} and show that, with the
previously introduced material, the $\mathcal{G}^{\mathcal{R}}$ wavefront of a
generalized function is defined exactly like the $\mathrm{C}^{\infty}$
wavefront of a distribution. First, as $\mathcal{G}^{\mathcal{R}}$ is a
subsheaf of $\mathcal{G}$, the following definition makes sense:

\begin{definition}
\label{RGDF-suppsing}Let $u$ be in $\mathcal{G}\left(  \Omega\right)  $. The
singular $\mathcal{G}^{\mathcal{R}}$-support of $u$ is the set
\[
\operatorname*{singsupp}\nolimits_{\mathcal{R}}u=\Omega\,\backslash\left\{
x\in\Omega\mathbb{\,}\left\vert \mathbb{\,}\exists V\in\mathcal{V}%
_{x},\;\right.  u\in\mathcal{G}^{\mathcal{R}}\left(  V\right)  \right\}
\text{.}%
\]

\end{definition}

\begin{proposition}
\label{RGDF-SpreS}$\mathcal{G}_{\mathcal{S}_{\ast}}:\Omega\rightarrow
\mathcal{G}_{\mathcal{S}_{\ast}}\left(  \Omega\right)  $ is a pre-sheaf: It
allows restrictions.
\end{proposition}

The \textbf{proof} is similar to the part $\left(  b\right)  $ of the one of
proposition \ref{GlimPropFRegS}.

\begin{notation}
\label{NotMicLoc}For $\left(  x,\xi\right)  \in\Omega\times\mathbb{R}%
^{d}\,\backslash\left\{  0\right\}  $ ($\Omega$ open subset of $\mathbb{R}%
^{d}$), we shall denote by:\newline$\left(  i\right)  ~\mathcal{V}_{x}$ (resp.
$\mathcal{V}_{\xi}^{\Gamma}$), the set of all open neighborhoods (resp. open
convex conic neighborhoods) of $x$ (resp. $\xi$),\newline$\left(  ii\right)
~\mathcal{D}_{x}\left(  \Omega\right)  $, the set of elements $\mathcal{D}%
\left(  \Omega\right)  $ non vanishing at $x$.
\end{notation}

For $\Gamma\in\mathcal{V}_{\xi}^{\Gamma}$, we say that $\widehat{u}%
\in\mathcal{G}_{\mathcal{S}_{\ast}}^{\mathcal{R}}\left(  \Gamma\right)  $ if
$u_{\left\vert \Gamma\right.  }\in\mathcal{G}_{\mathcal{S}_{\ast}%
}^{\mathcal{R}}\left(  \Gamma\right)  $. Let us fix a regular subset
$\mathcal{R}$ of $\mathbb{R}_{+}^{\mathbb{N}}$ and set, for $u\in
\mathcal{G}_{C}\left(  \mathbb{R}^{d}\right)  $,
\[
O^{\mathcal{R}}\left(  u\right)  =\left\{  \xi\in\mathbb{R}^{d}\,\backslash
\left\{  0\right\}  \,\left\vert \,\exists\Gamma\in\mathcal{V}_{\xi}^{\Gamma
}\;\ \right.  \widehat{u}\in\mathcal{G}_{\mathcal{S}_{\ast}}^{\mathcal{R}%
}\left(  \Gamma\right)  \right\}  \;\;\;\;\;\;\;\Sigma^{\mathcal{R}}\left(
u\right)  =\left(  \mathbb{R}^{d}\,\backslash\left\{  0\right\}  \right)
\backslash O^{\mathcal{R}}(u).
\]

\begin{lemma}
\label{RGDFLmHo811}For $u\in\mathcal{G}_{C}\left(  \mathbb{R}^{d}\right)  $
and $\varphi\in D\left(  \mathbb{R}^{d}\right)  ,\;O^{\mathcal{R}}\left(
u\right)  \subset O^{\mathcal{R}}\left(  \varphi u\right)  $ (or,
equivalently, $\Sigma^{\mathcal{R}}\left(  \varphi u\right)  \subset
\Sigma^{\mathcal{R}}\left(  u\right)  $).
\end{lemma}

\begin{proof}
Let $\left(  u_{\varepsilon}\right)  _{\varepsilon}\in\mathcal{X}\left(
\mathbb{R}^{d}\right)  $ be a representative of $u$ with $\operatorname*{supp}%
u_{\varepsilon}$ included in the same compact set, for all $\varepsilon$ in
$\left(  0,1\right]  $. We have
\[
\widehat{\varphi u}_{\varepsilon}(y)=\widehat{\varphi}\ast\widehat
{u}_{\varepsilon}(y)=\int\widehat{\varphi}\left(  \eta\right)  \widehat
{u}_{\varepsilon}\left(  y-\eta\right)  \,\mathrm{d}\eta.
\]
Let $\xi$ be in $O^{\mathcal{R}}\left(  u\right)  $ and $\Gamma\in
\mathcal{V}_{\xi}^{\Gamma}$ such that $\widehat{u}\in\mathcal{G}%
_{\mathcal{S}_{\ast}}^{\mathcal{R}}\left(  \Gamma\right)  $. There exists an
open conic neighborhood $\Gamma_{1}\subset\Gamma$ of $\xi$ and a real number
$c\in\left(  0,1\right)  $ such that, for all $\left(  y,\eta\right)  $ with
$y\in\Gamma_{1}$ and $\left\vert \eta\right\vert \leq c\left\vert y\right\vert
$, $y-\eta\in\Gamma$. Then
\begin{align*}
\widehat{\varphi u}_{\varepsilon}(y)  &  =\int_{\left\vert \eta\right\vert
\leq c\left\vert y\right\vert }\widehat{\varphi}\left(  \eta\right)
\widehat{u}_{\varepsilon}\left(  y-\eta\right)  \,\mathrm{d}\eta
+\int_{\left\vert \eta\right\vert >c\left\vert y\right\vert }\widehat{\varphi
}\left(  \eta\right)  \widehat{u}_{\varepsilon}\left(  y-\eta\right)
\,\mathrm{d}\eta\\
&  =\underset{v_{1,\varepsilon}(y)}{\underbrace{\int_{\left\vert
\eta\right\vert \leq c\left\vert y\right\vert }\widehat{\varphi}\left(
\eta\right)  \widehat{u}_{\varepsilon}\left(  y-\eta\right)  \,\mathrm{d}\eta
}}+\underset{v_{2,\varepsilon}(y)}{\underbrace{\int_{\left\vert y-\eta
\right\vert >c\left\vert y\right\vert }\widehat{\varphi}\left(  y-\eta\right)
\widehat{u}_{\varepsilon}\left(  \eta\right)  \,\mathrm{d}\eta}}.
\end{align*}
In order to estimate $v_{1,\varepsilon}$, let us remark that $\widehat{u}%
\in\mathcal{G}_{\mathcal{S}_{\ast}}^{\mathcal{R}}\left(  \Gamma\right)  $.
There exists a sequence $N\in\mathcal{R}$ such that, for all $q\in\mathbb{N}$,
there exists a constant $\ C_{1}>0$ with
\[
\forall\left(  y,\eta\right)  \in\Gamma_{1}\times\mathbb{R}^{d}\text{ with
}\left\vert \eta\right\vert \leq c\left\vert y\right\vert \;\;\;\;\left\vert
\widehat{u}_{\varepsilon}\left(  y-\eta\right)  \right\vert \leq
C_{1}\,\varepsilon^{-N\left(  q\right)  }\left(  1+\left\vert y-\eta
\right\vert \right)  ^{-q},
\]
for $\varepsilon$ small enough.

As, for $\left\vert \eta\right\vert \leq c\left\vert y\right\vert $, we have
$\left\vert y-\eta\right\vert \geq\left\vert \left\vert y\right\vert
-\left\vert \eta\right\vert \right\vert \geq\left\vert y\right\vert \left(
1-c\right)  $, it follows that%
\[
\forall\left(  y,\eta\right)  \in\Gamma_{1}\times\mathbb{R}^{d}\text{ with
}\left\vert \eta\right\vert \leq c\left\vert y\right\vert \;\;\;\;\left\vert
\widehat{u}_{\varepsilon}\left(  y-\eta\right)  \right\vert \leq
C_{1}\,\varepsilon^{-N\left(  q\right)  }\left(  1+\left\vert y\right\vert
\left(  1-c\right)  \right)  ^{-q}.
\]
Since $\widehat{\varphi}$ is rapidly decreasing, we get the existence of a
constant $C_{2}>0$ such that.
\[
\forall\eta\in\mathbb{R}^{d}\ \ \ \ \ \widehat{\varphi}\left(  \eta\right)
\leq C_{2}\left(  1+\left\vert \eta\right\vert \right)  ^{-d-1}.
\]
Replacing in the definition of $\left\vert v_{1,\varepsilon}(y)\right\vert $,
we get the existence of a constant $C_{3}>0$ such that
\[
\forall y\in\Gamma_{1}\ \ \ \ \left.  \left(  1+\left\vert y\right\vert
\right)  ^{q}\left\vert v_{1,\varepsilon}(y)\right\vert \right.  \leq
C_{3}\,\varepsilon^{-N\left(  q\right)  }\int\left(  \frac{1+\left\vert
y\right\vert }{\left(  1+\left\vert y\right\vert \left(  1-c\right)  \right)
}\right)  ^{q}\frac{1}{\left(  1+\left\vert \eta\right\vert \right)  ^{d+1}%
}\,\mathrm{d}\eta.
\]
The function $t\mapsto\left(  1+t\right)  /\left(  1+t\left(  1-c\right)
\right)  $ is bounded on $\mathbb{R}_{+}$.\ It follows that the integral in
the previous inequality converges, we finally get a constant $C_{4}>0$ such
that
\begin{equation}
\forall y\in\Gamma_{1}\ \ \ \ \ \left\vert v_{1,\varepsilon}(y)\right\vert
\ \leq C_{4}\,\varepsilon^{-N\left(  q\right)  }\left(  1+\left\vert
y\right\vert \right)  ^{-q}. \label{GregV1W}%
\end{equation}
For $v_{2,\varepsilon}$, note that $\left(  u_{\varepsilon}\right)
_{\varepsilon}\in\mathcal{X}_{\mathcal{S}_{\ast}}\left(  \mathbb{R}%
^{d}\right)  $.\ Therefore, there exist $M>0$ and $C_{5}>0$ such that
$\left\vert \widehat{u}_{\varepsilon}\left(  \eta\right)  \right\vert \leq
C_{5}\,\varepsilon^{-M}\left(  1+\left\vert \eta\right\vert \right)  ^{-d-1}$
for $\varepsilon$ small enough. As $\widehat{\varphi}\in S\left(
\mathbb{R}^{d}\right)  $, there exists $C_{6}>0$ such that
\[
\forall\left(  y,\eta\right)  \in\Gamma_{1}\times\mathbb{R}^{d}\text{ with
}\left\vert y-\eta\right\vert \geq c\left\vert y\right\vert \;\ \ \ \left\vert
\widehat{\varphi}\left(  y-\eta\right)  \right\vert \leq C_{6}\left(
1+\left\vert y-\eta\right\vert \right)  ^{-q}\leq C_{6}\left(  1+c\left\vert
y\right\vert \right)  ^{-q}.
\]
Then $\left\vert \widehat{\varphi}\left(  y-\eta\right)  \right\vert
=\mathrm{O}\left(  \left(  1+\left\vert y\right\vert \right)  ^{-q}\right)  $
as $y\rightarrow+\infty$. Thus, there exists a constant $C_{7}>0$ such that
\begin{equation}
\forall y\in\Gamma_{1}\ \ \ \ \ \ \left\vert v_{2,\varepsilon}(y)\right\vert
\leq C_{7}\,\varepsilon^{-M}\left(  1+\left\vert y\right\vert \right)
^{-q},\ \ \text{for }\varepsilon\text{ small enough.} \label{GregV2W}%
\end{equation}
From (\ref{GregV1W}) and (\ref{GregV2W}), we gets that, for all $q\in
\mathbb{N}$, there exists a constant $C>0$ (depending on $q$) such that
\[
\forall y\in\Gamma_{1}\ \ \ \ \ \ \left\vert \widehat{\varphi u}_{\varepsilon
}(y)\right\vert \leq C\,\varepsilon^{-\left(  N\left(  q\right)  +M\right)
}\left(  1+\left\vert y\right\vert \right)  ^{-q}.
\]
Since $\mathcal{R}$ is overstable by translation, there exists a sequence
$N^{\prime}\left(  \cdot\right)  \in\mathcal{R}$ such that $N\left(
\cdot\right)  +M\preceq N^{\prime}\left(  \cdot\right)  $ and $\mu
_{q,0}\left(  \widehat{\varphi u}_{\varepsilon}\right)  =\mathrm{O}\left(
\,\varepsilon^{-N^{\prime}\left(  q\right)  }\right)  $ as $\varepsilon
\rightarrow0$. Finally $\widehat{\varphi u}=\left[  \left(  \widehat{\varphi
u}_{\varepsilon}\right)  _{\varepsilon}\right]  _{\mathcal{G}_{\mathcal{S}%
_{\ast}}}\in\mathcal{G}_{\mathcal{S}_{\ast}}^{\mathcal{R}}\left(  \Gamma
_{1}\right)  $ and $\xi\in O^{\mathcal{R}}\left(  \varphi u\right)  $.
\end{proof}

\begin{definition}
\label{RGDF-DefMR}An element $u\in\mathcal{G}\left(  \Omega\right)  $ is said
to be $\mathcal{R}$ microregular on $\left(  x,\xi\right)  \in\mathbb{R}%
^{d}\times(\mathbb{R}^{d}\,\backslash\left\{  0\right\}  )$ if there exist
$\varphi\in\mathcal{D}_{x}\left(  \Omega\right)  $ and $\Gamma\in
\mathcal{V}_{\xi}^{\Gamma}$, such that $\widehat{\varphi u}\in\mathcal{G}%
_{\mathcal{S}_{\ast}}^{\mathcal{R}}\left(  \Gamma\right)  $.
\end{definition}

We set, for $u\in\mathcal{G}\left(  \Omega\right)  $ and $x\in\Omega$,%
\begin{align*}
O_{x}^{\mathcal{R}}\left(  u\right)   &  =\cup_{\varphi\in\mathcal{D}_{x}%
}O^{\mathcal{R}}\left(  \varphi u\right)  =\left\{  \xi\in(\mathbb{R}%
^{d}\,\backslash\left\{  0\right\}  )\,\left\vert \,u\text{ is microregular on
}\left(  x,\xi\right)  \right.  \right\}  ,\\
\Sigma_{x}^{\mathcal{R}}\left(  u\right)   &  =\cap_{\varphi\in\mathcal{D}%
_{x}}\Sigma^{\mathcal{R}}\left(  \varphi u\right)  =\left(  \mathbb{R}%
^{d}\,\backslash\left\{  0\right\}  \right)  \backslash O_{x}^{\mathcal{R}%
}\left(  u\right)  .
\end{align*}

\begin{definition}
\label{RGDF-WFHor}For $u\in\mathcal{G}\left(  \Omega\right)  $ the set
\[
WF_{\mathcal{R}}(u)=\left\{  \left(  x,\xi\right)  \in\mathbb{R}^{d}%
\times\mathbb{R}^{d}\backslash\left\{  0\right\}  \,\left\vert \,\xi\in
\Sigma_{x}^{\mathcal{R}}\left(  u\right)  \right.  \right\}
\]
is called the $\mathcal{R}$-wavefront of $u$.
\end{definition}

\begin{proposition}
\label{RGDFWF-SUPPS}For $u\in\mathcal{G}\left(  \Omega\right)  $, the
projection on the first component of $WF_{\mathcal{R}}(u)$ is equal to
$\operatorname*{singsupp}\nolimits_{\mathcal{R}}u$.
\end{proposition}

The \textbf{proof} of this proposition follows the same lines as the one for
the $\mathrm{C}^{\infty}$-wavefront of a distribution.\ First, for
$u\in\mathcal{G}\left(  \Omega\right)  $ and $\varphi\in\mathcal{D}\left(
\Omega\right)  $,\ $\varphi u$, which is \textit{a priori} in $\mathcal{G}%
_{C}\left(  \Omega\right)  $, can be straightforwardly considered as an
element of $\mathcal{G}_{C}\left(  \mathbb{R}^{d}\right)  $. As $\mathcal{G}%
_{C}\left(  \mathbb{R}^{d}\right)  $ is included in $\mathcal{G}_{\mathcal{S}%
}\left(  \mathbb{R}^{d}\right)  $ (see proposition \ref{GregEmbGCRGSRU}), the
Fourier transform of $\varphi u$ can be defined. (In the distributional case,
that is $u\in\mathcal{D}^{\prime}\left(  \Omega\right)  $, $\varphi u$ is
identified to an element of $\mathcal{E}^{\prime}\left(  \mathbb{R}%
^{d}\right)  $.) From this, we can follow the arguments of \cite{HorPDOT1}
(pages 253/254) for the $\mathrm{C}^{\infty}$-wavefront, which use mainly the
compactness of the sphere $S^{d-1}$ and lemma \ref{RGDFLmHo811}, which holds
in both cases: See lemma 8.1.1. in \cite{HorPDOT1} for the distributional case.

\begin{example}
\label{GregGinfWF}Taking $\mathcal{R}=\mathcal{B}$, the set of bounded
sequences, we recover the $\mathcal{G}^{\infty}$-wavefront, which has here a
definition independent of representatives.
\end{example}

\begin{example}
\label{GregG1WF}Taking $\mathcal{R}=\mathcal{R}_{1}$, we get a wavefront which
\textquotedblleft contains\textquotedblright\ the distributional microlocal
singularities of a generalized function, since $\mathcal{D}^{\prime}\left(
\cdot\right)  $ is embedded in $\mathcal{G}^{(1)}\left(  \cdot\right)  $.
\end{example}

In \cite{JAM3}, it is shown that the analogon of this lemma holds for the
analytic singularities of a generalized function, giving rise to the
corresponding wavefront set and the projection property of proposition
\ref{RGDFWF-SUPPS}. Our future aim is to apply this theory to the propagation
of singularities through integral generalized operators \cite{DELGST}. We also
refer the reader to \cite{GarHorm}, \cite{HorDeH}, \cite{HorKun},
\cite{HorObPil}, \cite{NePiSc} and the literature therein for other
presentations of the $\mathcal{G}^{\infty}$-wavefront (which is a particular
case of $\mathcal{R}$ wavefront).

\bigskip

\end{document}